\documentclass[11pt,english,a4paper,twoside]{article}
\usepackage{hyperref}
\usepackage{amsthm,amsfonts,amssymb,a4wide}
\usepackage[pagestyles]{titlesec}
\usepackage{lipsum}
\usepackage{babel}
\usepackage[utf8]{inputenc}
\usepackage[T1]{fontenc}
\usepackage[many]{tcolorbox}
\tcbuselibrary{listings,theorems,breakable,hooks,xparse,skins,fitting}
\usepackage{makeidx}
\usepackage{showidx}
\usepackage{wasysym}
\usepackage{calligra}
\usepackage{cancel}
\usepackage{pst-plot,pstricks-add,pstricks,pst-coil,pst-node,pst-func,pst-math,pst-xkey,pst-3dplot,pst-lens,pst-solides3d}
\usepackage{transparent}
\usepackage{graphicx,textcomp,pifont,shapepar}
\usepackage{color}
\usepackage{polynom}
\usepackage[tiling]{pst-fill}
\usepackage{array}
\usepackage{rotating}
\usepackage{mdframed}
\usepackage{multicol}
\usepackage{tikz}
\usetikzlibrary{shapes,arrows}
\usepackage{wallpaper}
\usepackage{watermark}
\usepackage{fancyhdr}
\usepackage{transparent}
\usepackage{enumerate}
\usepackage[cam,center]{crop}
\usepackage{colortbl}
\usepackage{xpatch}
\usepackage{appendix}
\usepackage{cleveref}
\usepackage{tcolorbox}
\usepackage{changes}
\usepackage{soul}
\usepackage{cancel}

\makeindex
\usepackage[T1]{fontenc}

\AtBeginDocument{
\let\div\undefined\DeclareMathOperator{\div}{div}
}

\DeclareMathOperator{\grad}{grad}

\newcommand\numberthis{\addtocounter{equation}{1}\tag{\theequation}}

\newcommand{\e}{\text{e}}

\newcommand{\p}{\pi_{\mathcal{N}}} 

\newcommand{\f}{\text{d}} 

\newcommand{\si}{\psi_{\nu,\tau}} 

\newcommand{\B}{B_{\tau}^+}

\newcommand{\N}{\mathbb{N}} 

\newcommand{\R}{\mathbb{R}}

\newcommand{\U}{\mathcal{U}} 

\newcommand{\h}{\mathcal{H}} 

\newcommand{\s}{\text{supp}} 

\newtheorem{thm}{Theorem}[section]

\newtheorem{defn}[thm]{Definition}

\newtheorem{lem}[thm]{Lemma}

\newtheorem{remark}[thm]{Remark}

\titleformat{\chapter}[display]
        {\color{black} \Huge\bfseries}{\chaptertitlename\ \thechapter}{0pt}{\Huge\bfseries}

\titleformat{\section}[block]{\color{black}\Large\bfseries}{\thesection\quad}{0pt}{\Large}

\titleformat{\subsection}[block]
        {\color{black}\large\bfseries}{\thesubsection}{10pt}{\large}

\def\XXint#1#2#3{{\setbox0=\hbox{$#1{#2#3}{%
\int}$ }
\vcenter{\hbox{$#2#3$ }}\kern-.6\wd0}}

\begin{document}
\title{\Large\textbf{Derivation of a boundary monotonicity inequality for variationally biharmonic maps}}
\author{Serdar Altuntas\footnote{Fakultät für Mathematik, Universität Duisburg-Essen, 45117 Essen, Germany}}
\date{ }
\maketitle
\begin{abstract}
We derive a boundary monotonicity formula for a class of biharmonic maps with Dirichlet boundary conditions. A monotonicity formula is crucial in the theory of partial regularity in super-critical dimensions. As a consequence of such a boundary monotonicity formula, one is able to show partial regularity for variationally biharmonic maps and full boundary regularity for minimizing biharmonic maps.
\end{abstract}

\section{Introduction}
Over the last decades it has turned out that a monotonicity formula is necessary in super-critical dimensions to show partial regularity. Before the study of weakly biharmonic maps has begun, one has considered weakly harmonic maps. Let $\mathcal{M}$ be a smooth Riemannian manifold of dimension $m\in \N$ with or without boundary and $\mathcal{N}\subset \R^n$ be a compact Riemannian manifold without boundary. We call a map $u\in W^{1,2}(\mathcal{M},\mathcal{N})$ weakly harmonic iff it is a critical point of the so-called Dirichlet-energy
\begin{align*}
E_1(u)=\int_{\mathcal{M}}\vert Du\vert^2\f \mu_{\mathcal{M}},\numberthis \label{1}
\end{align*}
for variations of the form $u_t=\p(u+tV)$ for $V\in C^{\infty}_0(\mathcal{M},\R^n)$. Here, $\p$ denotes the nearest point projection. Critical points of $E_1$ satisfy a nonlinear system of second order equations
\begin{align*}
\Delta u=\text{tr}(A\circ u)(Du\otimes Du)\numberthis \label{1b}
\end{align*}
in the sense of distribution with a critically nonlinear right-hand side where $\text{tr}(A)$ denotes the trace of the second fundamental form of $\mathcal{N}$. There are several regularity results of weakly harmonic maps. In 1948 C.B. Morrey \cite{22} showed that every minimizing map $u\in W^{1,2}(\mathcal{M},\mathcal{N})$ is $C^{\infty}$ for a manifold of dimension $\dim{\mathcal{M}}=m\leq2$. For $m=2$, F. Héléin \cite{23} proved that any weakly harmonic map $u\in W^{1,2}(\mathcal{M},\mathcal{N})$ is smooth inside $\mathcal{M}$. The right-hand side is a priori just in $L^1(\mathcal{M},\mathcal{N})$. Therefore, the information from \eqref{1b} is not enough to get some regularity results in dimensions $m>2$. A counter-example of T. Riviére \cite{24} illustrates this fact. In 1995 he constructed an everywhere discontinuous weakly harmonic map. Therefore, one has to consider stationary harmonic maps which are weakly harmonic and in addition critical points of $E_1$ for inner variations. A useful property of stationary harmonic maps is that they fulfil an energy monotonicity formula which is crucial to show partial regularity in super-critical dimensions. The first result of partial regularity for stationary harmonic maps in arbitrary compact manifolds was shown by Bethuel \cite{25} which is a generalisation of Evans work in \cite{6} where he considered maps from a subset of the Euclidean space into the unit sphere $\mathcal{N}=\mathbb{S}^{n-1}$. Another class of harmonic maps are energy minimizing harmonic maps. We call $u\in W^{1,2}(\mathcal{M},\mathcal{N})$ a minimizing harmonic map if $E_1(u)\leq E_1(v)$ for all $v\in W^{1,2}(\mathcal{M},\mathcal{N})$ such that $u-v\in W_0^{1,2}(\mathcal{M},\mathcal{N})$. R. Schoen und K. Uhlenbeck \cite{26,27} established interior partial regularity and boundary regularity for minimizing harmonic maps. An analogy to weakly harmonic maps are (extrinsically\footnote{One distinguishes between extrinsically and intrisically biharmonic maps. We say that a map is intrinsically biharmonic iff it is a critical point of $\mathcal{E}(u)=\int_{\mathcal{M}}\vert \nabla Du\vert^2\f\mu_{\mathcal{M}}$. The energy $\mathcal{E}$ does not depend on the embedding $\mathcal{N}\hookrightarrow \R^n$ while $E_2$ does. Therefore, the distinction extrinsically and intrinsically.}) weakly biharmonic maps which are critical points of the so-called bienergy or Hessian energy
\begin{align*}
E_2(u)=\int_{\mathcal{M}}\vert \Delta u\vert^2\f \mu_{\mathcal{M}}.\numberthis \label{1c}
\end{align*} 
They were firstly studied by S.-Y. A. Chang, L. Wang and P. C. Yang in \cite{4} in domains of dimension greater than or equal four into spheres. Again, a monotonicity formula for stationary biharmonic maps in super-critical dimensions was crucial to show interior partial regularity. However, they derived this monotonicity formula only for sufficiently regular maps. G. Angelsberg \cite{1} gave a rigorous proof of this monotonicity formula for stationary biharmonic maps $u\in W^{2,2}(B_r,\mathcal{N})$. A monotonicity formula for intrinsically stationary biharmonic maps was derived by R. Moser \cite{30}. In the case of minimizing maps, M.-C. Hong and C. Wang \cite{28} showed that any minimizing biharmonic map for $\mathcal{N}=\mathbb{S}^{n-1}$ is smooth off a singular set $\Sigma$ whose Hausdorff dimension is at most $m-5$, where $m\in \N_{\geq 5}$. C. Scheven \cite{3} showed that for an arbitrary target manifold $\mathcal{N}$ the singular set of a minimizing biharmonic map has Hausdorff dimension at most $m-5$. A boundary regularity theory for stationary biharmonic maps was initiated by H. Gong, T. Lamm and C. Wang in \cite{7}. They derived a boundary monotonicity inequality for biharmonic maps of class $ W^{4,2}(\Omega,\mathcal{N})$, where $\Omega=B_R^+(a)$ is a half ball with Euclidian metric. Both assumptions are not natural. The first assumption '$W^{4,2}$' trivializes the regularity problem. The second assumption '$\Omega=B_R^+(a)$' excludes curved parts of the boundary. Therefore, a flattening of the boundary will change the bienergy functional $E_2$ by lower order terms. Furthermore, K. E. Mazowiecka \cite{29} proved recently in her dissertation that minimizing biharmonic maps are smooth in a full neighborhood of the boundary under the assumption that there exists a boundary monotonicity formula. However, the proof of the boundary monotonicity inequality is missing and this turns out to be technically very demanding. We derive in \Cref{sec3} a boundary monotonicity inequality for a class of biharmonic maps in the function space $W^{2,2}(B_R^+,\mathcal{N})$ and close this gap in Mazowiecka's dissertation. In this sense, we provide the last missing ingredient for the proof of the full boundary regularity of minimizing biharmonic maps. We also include the case of a curved boundary. We proceed as in \cite[Theorem 2]{17}, i.e. we consider variations of the form $u_t=\p(u\circ \varphi_t-g\circ \varphi_t+g)$ and use the methods in \cite{1}. Since we allow slightly more general variations than in the case of stationary biharmonic maps, we call our maps \textit{variationally biharmonic maps} similarly to \cite{18}. For the derivation of the boundary monotonicity inequality \eqref{57} we need at first a differential equation which we derive in \Cref{sec2}.\\
Now, we introduce our setting and give some definitions: Let $\Omega\subset \R^m$ be a bounded domain with smooth boundary of dimension $m\geq5$ equipped with a smooth Riemannian metric $\gamma$ and $\mathcal{N}$ be a smooth, compact Riemannian manifold without boundary which is isometrically embedded in Euclidean space $\R^n$. For $$u\in W^{2,2}(\Omega,\mathcal{N}):=\left\lbrace u\in W^{2,2}(\Omega,\R^n): u(x)\in \mathcal{N} \text{ for a.e. } x\in \Omega\right\rbrace$$ satisfying Dirichlet boundary conditions $$\left(u,Du\right)\vert_{\partial \Omega}=\left(g,Dg\right)\vert_{\partial \Omega}$$ in the sense of trace for given boundary data $g\in C^3(\Omega,\mathcal{N})$ the so-called extrinsic bienergy functional is defined as
\begin{align*}
E(u)=\int_{\Omega}\vert \Delta_{\gamma}u\vert^2\f\mu_{\gamma}.\numberthis \label{2}
\end{align*}
Here, $\Delta_{\gamma}:=\gamma^{ij}\left(\partial_i\partial_j-\Gamma_{ij}^k\partial_k\right)$ denotes the \textit{Laplace-Beltrami-operator} and $\mu_{\gamma}:=\mathcal{L}^m\llcorner\sqrt{\gamma}$ stands for the Riemannian measure on $\Omega$, where $\sqrt{\gamma}:=\sqrt{\det{(\gamma_{ij})}}$ and $\Gamma_{ij}^k:=\frac{1}{2}\gamma^{kl}(\partial_i\gamma_{jl}-\partial_l\gamma_{ij}+\partial_j\gamma_{il})$ are the Christoffel-symbols of the second kind.

The \textit{Riemannian gradient} $\grad_{\gamma} f(x)$ of $f\in C^1(\Omega,\R)$ is defined by $\gamma\left(\grad_{\gamma}f(x),X\right)=X(f)$ for all $x\in \Omega$ and every vector field $X=X^i\partial_i\in C^1(\Omega,\R^n)$. In coordinates we have $\grad_{\gamma}f(x)=\gamma^{ij}(x)\partial_if(x)\partial_j$. The \textit{Riemannian divergence} $\div_{\gamma}$ of a vector field $X\in C^1(\Omega,\R^n)$ is defined as the trace of the map $Y\mapsto \nabla_{Y}X$, where $\nabla$ denotes the covariant derivative. In coordinates, $\div_{\gamma}{X}=\dfrac{1}{\sqrt{\gamma}}\partial_k\left( \sqrt{\gamma}X^k\right)$.

For $\delta>0$, let  $V_{\delta}$ be a neighborhood of $\mathcal{N}$, which is given by $V_{\delta}:=\left\lbrace p\in \R^n: \text{dist}(p,\mathcal{N})<\delta\right\rbrace $. Since $\mathcal{N}$ is smooth and compact, there are sufficiently small $\delta>0$, so that for all $p\in V_{\delta}$ a unique point $\p(p)\in\mathcal{N}$ with $\vert p-\p(p)\vert =\text{dist}(p,\mathcal{N})$ exists. The map $\p:V_{\delta}\rightarrow \mathcal{N}$ is called \textit{nearest point projection}. The total derivative of $\p$ in $p\in \mathcal{N}$ is the \textit{orthogonal projection} onto the tangential space in $p$, i.e. $D\p: \R^n\rightarrow T_p\mathcal{N}$. For more details see for example Moser \cite[chapter 3]{15} or \cite[chapter 2.12.3]{12}.

A map $u\in W^{2,2}(\Omega,\mathcal{N})$ is said to be \textit{weakly biharmonic} if and only if it satisfies
\begin{align*}
\dfrac{\f}{\f t}\bigg\vert_{t=0} E(u_t)=0\numberthis \label{2b}
\end{align*}
for all variations of the form $u_t=\p(u+t\psi)$ with $\psi\in C^{\infty}_0(\Omega,\R^n)$. A weakly biharmonic map $u\in W^{2,2}(\Omega,\mathcal{N})$ is called \textit{stationary biharmonic} if it satisfies \eqref{2b} additionally for variations of the form $u_t(x)=u(x+t\xi(x))$ with $\xi\in C^{\infty}_0(\Omega,\R^m)$. We say that $u\in W^{2,2}(\Omega,\mathcal{N})$ is a \textit{minimizing biharmonic} map if and only if $E(u)\leq E(v)$ for all $v\in W^{2,2}(\Omega,\mathcal{N})$ with $u-v\in W_0^{2,2}(\Omega,\mathcal{N})$. Clearly, minimizing biharmonic maps are stationary biharmonic. We give another class of biharmonic maps in the following 
\begin{defn}\label{def}
We name a map $u\in W^{2,2}(\Omega,\mathcal{N})$ \textbf{variationally biharmonic with respect to the Dirichlet boundary conditions} $\left(u,Du\right)\vert_{\partial \Omega}=\left(g,Dg\right)\vert_{\partial \Omega}$ if it is weakly biharmonic and satisfies \eqref{2b} for variations of the form $u_t(x)=\p(u\circ \varphi_t-g\circ \varphi_t+g)$. Here, $\varphi_t$ is a $C^{\infty}$-family of diffeomorphisms from $\Omega$ into $\Omega$ that satisfy $\varphi(\partial\Omega)\subset \partial\Omega$, $\varphi_0=id_{\Omega}$, $\dfrac{\f}{\f t}\bigg\vert_{t=0}\varphi_t(x)=\xi(x)\in C^{\infty}(\Omega,\R^m)$.
\end{defn}

Throughout, we use the following notations
\begin{align*}
\int_{\partial B_r^+\setminus \partial B_{\rho}^+}f:=\int_{\partial B_r^+}f-\int_{\partial B_{\rho}^+}f,\quad & \int_{S_{r}^+\cup S_{\rho}^+}f:=\int_{S_{r}^+}f+\int_{S_{\rho}^+}f.
\end{align*}

Now, we state our main result:
\begin{thm}[\textbf{Boundary monotonicity inequality}]\label{thm}
For $m\in \N_{\geq5}$, let $u\in W^{2,2}(B^+,\mathcal{N})$ be a variationally biharmonic map from the half-ball $B^+:=B_R^+(a):=B_R(a)\cap \R^{m-1}\times[0,\infty)$ with center $a\in \R^{m-1}\times \left\lbrace 0\right\rbrace $ and radius $R>0$ to a Riemannian manifold $\mathcal{N}\subset \R^n$. Let $B^+$ be equipped with a general smooth Riemannian metric $\gamma$, where the metric satisfies $\gamma_{ij}(a)=\delta_{ij}$ for $1\leq i,j\leq n$ and for an ellipticity constant $G\geq1$ and a constant $H\geq 0$ the conditions
\begin{align*}
G^{-1}\vert \theta\vert^2\leq \sum_{i,j=1}^m\gamma_{ij}(x)\theta^i\theta^j\leq G\vert \theta\vert^2, \qquad \vert \gamma_{ij}(x)\vert\leq G, \qquad \vert \partial_k\gamma_{ij}(x)\vert,\vert \partial_l\partial_k\gamma_{ij}(x)\vert\leq H \numberthis \label{1}
\end{align*}
for all  $x\in B^+$, $\theta=(\theta^1,\ldots,\theta^m)\in \R^m$. Furthermore, we denote the curved and flat part of $\partial B^+$ by $S_R^+:=\partial B_R\cap \left\lbrace  x\in \R^m:x^m>0\right\rbrace$ and $T_R:=\partial B_R^+\cap \left\lbrace  x\in \R:x^m=0\right\rbrace$. Suppose that the Dirichlet boundary conditions $\left(u,Du\right)\vert_{T_R}=\left(g,Dg\right)\vert_{T_R}$ hold for given boundary data $g\in C^3(B^+,\mathcal{N})$. Then, there are constants $\chi=\chi(\mathcal{N},G,H,\Vert Dg\Vert_{C^2},\Vert u\Vert_{L^{\infty}(B^+_1)})\geq0$, $\mathsf{C}_1,\ldots,\mathsf{C}_7\geq0$, so that
 \begin{align*}
&\e^{\chi r}r^{4-m}\int_{B_r^+(a)}\vert \Delta_{\gamma}u\vert^2 \f \mu_{\gamma}-\e^{\chi\rho}\rho^{4-m}\int_{B_{\rho}^+(a)}\vert \Delta_{\gamma}u\vert^2 \f \mu_{\gamma}+\mathsf{C}_1 r\\
&\quad+\mathsf{C}_2\int_{\rho}^r\e^{\chi\tau}\tau^{5-m}\int_{B^+_{\tau}(a)}\vert D^2u\vert^2 \f\mu_{\gamma}\f\tau+ \mathsf{C}_3\int_{\rho}^r\e^{\chi\tau}\tau^{3-m}\int_{B^+_{\tau}(a)} \vert Du\vert^2 \f\mu_{\gamma}\f\tau\\
&\quad +\mathsf{C}_4\int_{ B^+_{r}(a)\setminus B^+_{\rho}(a)}\e^{\chi\vert x-a\vert}\dfrac{\vert D^2 u\vert^2  }{\vert x-a\vert^{m-5}}\f\mu_{\gamma}+\mathsf{C}_5\int_{ B^+_{r}(a)\setminus B^+_{\rho}(a)}\e^{\chi\vert x-a\vert}\dfrac{\vert Du\vert^2  }{\vert x-a\vert^{m-3}}\f\mu_{\gamma}\\
&\quad+\mathsf{C}_6\int_{S_{r}^+(a)\cup S_{\rho}^+(a)}\e^{\chi\vert x-a\vert}\dfrac{\vert D^2u\vert^2  }{\vert x-a\vert^{m-6}} \sqrt{\gamma}\f\h^{m-1}\\
&\quad+\mathsf{C}_7\int_{S_{r}^+(a)\cup S_{\rho}^+(a)}\e^{\chi\vert x-a\vert}\dfrac{\vert Du\vert^2}{\vert x-a\vert^{m-4}}\sqrt{\gamma} \f\h^{m-1}\\
&\geq 4\int_{B^+_r(a)\setminus B^+_{\rho}(a)}\e^{\chi\vert x-a\vert}\left( \dfrac{ (u_i+u_{ij}(x-a)^j)^2}{\vert x-a\vert^{m-2}}+\dfrac{(m-2)\vert Du\cdot (x-a)\vert^2}{\vert x-a\vert^{m}}\right)\f\mu_{\gamma}\label{57}\tag{M}\\
&\quad+2\int_{S^+_r(a)\setminus S^+_{\rho}(a)}\e^{\chi\vert x-a\vert}\left(-\dfrac{ u_{i} u_{ij}(x-a)^j}{\vert x-a\vert^{m-3}}+ 2\dfrac{ \vert Du\cdot (x-a)\vert^2}{\vert x-a\vert^{m-1}}-2\dfrac{ \vert Du\vert^2}{\vert x-a\vert^{m-3}}\right)\sqrt{\gamma}\f\h^{m-1}
\end{align*}
holds for a.e. $0<\rho<r<R$, where $\mathsf{C}_1,\ldots,\mathsf{C}_7$ depend on $m,\mathcal{N}, G, H$ and $\Vert Dg\Vert_{C^2}$ and  $\mathsf{C}_1,\ldots,\mathsf{C}_5$ additionally on $\Vert u\Vert_{L^{\infty}(B^+_1)}$. Moreover, $\chi$ and $\mathsf{C}_1$ to $\mathsf{C}_7$ vanish for $Dg\rightarrow0$ in $C^2$ and for constant metric $\gamma$.
\end{thm}

There are two known consequences of such a boundary monotonicity inequality similar to \eqref{57}. The first one was shown by H. Gong, T. Lamm and C. Wang \cite{7}. They obtained the following result: if $u$ is a stationary biharmonic map that satisfies a certain boundary monotonicity inequality, then there exists a closed subset $\Sigma\subset\bar{\Omega}$, with $\h^{m-4}(\Sigma)=0$, such that $u\in C^{\infty}(\bar{\Omega}\setminus \Sigma,\mathcal{N})$. The second one was established by K. Mazowiecka \cite{29}. She proved that every minimizing biharmonic map which satisfies a certain boundary inequality is smooth on a full neighborhood of the boundary $\partial\Omega$. In both proofs, $\Omega$ is a subset of $\R^m$, $m\geq5$, with Euclidean metric and $\left(u,\dfrac{\partial u}{\partial \nu}\right)\bigg\vert_{\partial\Omega}=\left(\phi,\dfrac{\partial \phi}{\partial \nu}\right)\bigg\vert_{\partial\Omega}$ for a given $\phi\in C^{\infty}(\Omega_{\delta},\mathcal{N})$ where $\Omega_{\delta}=\left\lbrace x\in \bar{\Omega}: \text{dist}(x,\partial\Omega)<\delta\right\rbrace $ for some $\delta>0$. Here, $\nu$ denotes the outer normal vector.

\section{Differential equation for variational biharmonic maps}\label{sec2}
The starting point for our derivation of the boundary monotonicity inequality \eqref{57} is the differential equation \eqref{7} in the following
\begin{lem}\label{lemma1}
Let $u\in W^{2,2}(\Omega,\mathcal{N})$ be a variational biharmonic with respect to the Dirichlet boundary conditions $\left(u,Du\right)\vert_{\partial \Omega}=\left(g,Dg\right)\vert_{\partial \Omega}$, then the following differential equation holds for all $\xi\in C^{\infty}(\Omega,\R^m)$ with $\xi\in T_x(\partial \Omega)$ for every $x\in \partial\Omega$:
\begin{align*}
&\int_{\Omega}\left(4\Delta_{\gamma}u\cdot D^2u\grad_{\gamma}\xi+2\Delta_{\gamma}u\cdot Du\Delta_{\gamma}\xi-\vert \Delta_{\gamma}u\vert^2\div_{\gamma}\xi\right)\f\mu_{\gamma}\label{7}\tag{D} \\
&=\int_{\Omega}\left(2\Delta_{\gamma}u\cdot \Delta_{\gamma}\left[D\p(u)\left(Dg\xi\right)\right]+2\Delta_{\gamma}u\cdot \partial_l\gamma^{ij}\xi^l\partial_i\partial_ju-2\Delta_{\gamma}u\cdot \partial_l\left(\gamma^{ij}\Gamma_{ij}^k\right)\xi^l\partial_ku\right)\f\mu_{\gamma}
\end{align*}
\end{lem}
Here, '$\partial_i$' denotes partial derivation with respect to $x^i$.
\begin{proof}
Let $\varphi_t$ be as in \Cref{def} with $\dfrac{\f}{\f t}\bigg\vert_{t=0}\varphi_t=\xi\in C^{\infty}(\Omega,\R^m)$. We consider the variation $\U_t(x):=u\circ \varphi_t(x)-g\circ \varphi_t(x)+g(x)$ in $\R^n$. For $x\in \partial\Omega$ it holds $\U_t(x)=g(x)$ and $\partial_l\U_t(x)=\partial_lg(x)$ for all $1\leq l\leq m$. So, $\U_t$ satisfies the boundary conditions, and it holds $\U_0=u$. Since $u\circ \varphi_t(x)\in \mathcal{N}$, the image of $\U_t(x)$ is for sufficiently small $\vert t\vert$ in a neighborhood of $\mathcal{N}$, i.e. in the domain of $\pi_{\mathcal{N}}$. Thus, we consider the variation $x\mapsto \p(\U_t(x))=:u_t(x)$ and therefore the following functional
\begin{align*}
\int_{\Omega}\vert \Delta_{\gamma} u_t(x)\vert^2\f\mu_{\gamma}.\numberthis \label{8}
\end{align*}
With the transformation $x\mapsto \varphi^{-1}_t(x)$ we get
\begin{align*}
\int_{\Omega}\vert \Delta_{\gamma} u_t(x)\vert^2\f\mu_{\gamma}&=\int_{\Omega}\vert \Delta_{\gamma}u_t\circ\varphi_t^{-1}(x)\vert^2 \det{D\varphi_t^{-1}(x)}\sqrt{\gamma\circ\varphi_t^{-1}(x)}\f\mathcal{L}^m(x)\\
&:=\int_{\Omega}f(t,x)\f\mathcal{L}^m.\numberthis \label{9}
\end{align*}
To derive the equation \eqref{7} we differentiate the functional \eqref{9} with respect to $t$ and evaluate the result at $t=0$. Since we consider variational biharmonic maps, it holds
\begin{align*}
0&=\dfrac{\f}{\f t}\bigg\vert_{t=0}\int_{\Omega}f(t,x)\f\mathcal{L}^m.\numberthis \label{9b}
\end{align*}
For the sake of clarity we omit the argument '$x$'. Now, it holds
\begin{align*}
\dfrac{\f}{\f t}\bigg\vert_{t=0}f(t,x)&=2 \Delta_{\gamma}u\cdot  \left(\dfrac{\f}{\f t}\bigg\vert_{t=0}\Delta_{\gamma}u_t\circ\varphi_t^{-1}\right) \sqrt{\gamma}+\vert \Delta_{\gamma}u\vert^2\cdot \left(\dfrac{\f}{\f t}\bigg\vert_{t=0} \det{D\varphi_t^{-1}} \right)\sqrt{\gamma}\\
&\quad+\vert \Delta_{\gamma}u\vert^2\cdot \left(\dfrac{\f}{\f t}\bigg\vert_{t=0}\sqrt{\gamma\circ\varphi_t^{-1}}\right)\numberthis\label{10}.
\end{align*}
In six steps we compute the following three terms, $
\dfrac{\f}{\f t}\bigg\vert_{t=0}\Delta_{\gamma}u_t\circ\varphi_t^{-1}$, 
$\displaystyle \dfrac{\f}{\f t}\bigg\vert_{t=0} \det{D\varphi_t^{-1}} $ and $\displaystyle \dfrac{\f}{\f t}\bigg\vert_{t=0}\sqrt{\gamma\circ\varphi_t^{-1}}$:\\

\textit{Step 1}. We have 
\begin{align*}
\Delta_{\gamma}u_t=\gamma^{ij}\left(\partial_i\partial_ju_t-\Gamma^k_{ij}\partial_ku_t \right)\numberthis\label{12}
\end{align*}
where 
\begin{align*}
\partial_ku_t=\partial_k\left[\p(\U_t)\right]=D\p(\U_t)(\partial_k\U_t)\numberthis\label{13}
\end{align*}
with
\begin{align*}
\partial_k\U_t=\left( \partial_l(u-g)\circ \varphi_t  \right)\partial_k\varphi_t^l+\partial_kg\numberthis\label{14}
\end{align*} 
and 
 \begin{align*}
 \partial_i\partial_ju_t&=\partial_i\partial_j\left[\p(\U_t)\right] = D^2\p(\U_t)(\partial_i\U_t,\partial_j\U_t)+D\p(\U_t)(\partial_i\partial_j\U_t)\numberthis\label{15}
 \end{align*}
 with
 \begin{align*}
 \partial_i\partial_j\U_t&=\left( \partial_k\partial_l(u-g)\circ \varphi_t  \right)\partial_i\varphi_t^l\partial_j\varphi_t^k+\left( \partial_k(u-g)\circ \varphi_t  \right)\partial_i\partial_j\varphi_t^k+\partial_i\partial_jg.\numberthis\label{16}
 \end{align*}
Now, we get from \eqref{12}
\begin{align*}
\Delta_{\gamma}u_t\circ\varphi_t^{-1}&=\gamma^{ij}\circ \varphi_t^{-1}\left(\partial_i\partial_ju_t\circ \varphi_t^{-1}-\Gamma^k_{ij}\circ \varphi_t^{-1}\partial_ku_t\circ \varphi_t^{-1} \right).\numberthis\label{17}
\end{align*}
Due to \eqref{13} and \eqref{15} we have
\begin{align*}
\partial_ku_t\circ \varphi_t^{-1}=D\p(\U_t\circ \varphi_t^{-1})(\partial_k\U_t\circ \varphi_t^{-1})\numberthis\label{18}
\end{align*}
and
 \begin{align*}
 \partial_i\partial_ju_t\circ \varphi_t^{-1}&= D^2\p(\U_t\circ \varphi_t^{-1})(\partial_i\U_t\circ \varphi_t^{-1},\partial_j\U_t\circ \varphi_t^{-1})\\
 &\quad+D\p(\U_t\circ \varphi_t^{-1})(\partial_i\partial_j\U_t\circ \varphi_t^{-1}).\numberthis\label{19}
 \end{align*} 
 
\textit{Step 2}. From \eqref{17} we obtain by using the product rule,
\begin{align*}
\dfrac{\f}{\f t}\bigg\vert_{t=0}\Delta_{\gamma}u_t\circ\varphi_t^{-1}&=\left(\dfrac{\f}{\f t}\bigg\vert_{t=0}\gamma^{ij}\circ \varphi_t^{-1}\right)\cdot\left(\partial_i\partial_ju-\Gamma_{ij}^k\partial_ku\right)+\gamma^{ij}\cdot\left(\dfrac{\f}{\f t}\bigg\vert_{t=0}\partial_{i}\partial_{j}u_t\circ \varphi_t^{-1}\right)\\
&\quad-\gamma^{ij}\cdot\left[ \left(\dfrac{\f}{\f t}\bigg\vert_{t=0}\Gamma^k_{ij}\circ \varphi_t^{-1}\right) \cdot \partial_ku+\Gamma^k_{ij}\cdot \left(\dfrac{\f}{\f t}\bigg\vert_{t=0}\partial_{k}u_t\circ \varphi_t^{-1} \right)\right].\numberthis\label{20}
\end{align*} 

\textit{Step 3}. Next, we compute $\dfrac{\f}{\f t}\bigg\vert_{t=0}\partial_{j }u_t\circ \varphi_t^{-1} $ and $ \dfrac{\f}{\f t}\bigg\vert_{t=0}\partial_{i}\partial_{j}u_t\circ \varphi_t^{-1}$. Due to \eqref{18} we have
\begin{align*}
\dfrac{\f}{\f t}\bigg\vert_{t=0}\partial_{k}u_t\circ \varphi_t^{-1} &=D^2\p(u)\left(\partial_ku,\dfrac{\f}{\f t}\bigg\vert_{t=0}\U_t\circ\varphi_t^{-1}\right)+D\p(u)\left(\dfrac{\f}{\f t}\bigg\vert_{t=0}\partial_k\U_t\circ\varphi_t^{-1}\right)\numberthis\label{21}.
\end{align*}
Since $\U_t\circ\varphi_t^{-1}=(u-g)+g\circ \varphi_t^{-1}$, it holds
\begin{align*}
\dfrac{\f}{\f t}\bigg\vert_{t=0}\U_t(\varphi_t^{-1})&=Dg\dfrac{\f}{\f t}\bigg\vert_{t=0}\varphi_t^{-1}=-Dg\xi,\numberthis\label{22}
\end{align*}
where we used in the last step that
\begin{align*}
\dfrac{\f}{\f t}\bigg\vert_{t=0}\varphi_t^{-1}=-\xi,\numberthis\label{24}
\end{align*}
which is a consequence of the chain rule and the fact $\dfrac{\f}{\f t}\varphi_t=\xi$. Equation \eqref{14} yields
\begin{align*}
\partial_k\U_t\circ\varphi_t^{-1}=\left( \partial_l(u-g)  \right)\partial_k\varphi_t^l\circ\varphi_t^{-1}+\partial_kg\circ\varphi_t^{-1}\numberthis\label{26}.
\end{align*} 
Consequently, we get with \eqref{24} the equation
\begin{align*}
\dfrac{\f}{\f t}\bigg\vert_{t=0}\partial_k\U_t\circ\varphi_t^{-1}&=\left( \partial_l(u-g) \right)\dfrac{\f}{\f t}\bigg\vert_{t=0}\left(\partial_k\varphi_t^l\circ\varphi_t^{-1}\right)+\dfrac{\f}{\f t}\bigg\vert_{t=0}\left(\partial_kg\circ\varphi_t^{-1}\right)\\
&=\partial_l(u-g) \dfrac{\f}{\f t}\bigg\vert_{t=0}\left(\partial_k\varphi_t^l\circ\varphi_t^{-1}\right)-\partial_l\partial_kg\xi^l\\
&=\partial_l(u-g)\partial_k\xi^l-\partial_l\partial_kg\xi^l=D(u-g)\partial_k\xi-\partial_kDg\xi\numberthis\label{29}.
\end{align*}
We put the equations \eqref{22} and \eqref{29} into \eqref{21}, and obtain
\begin{align*}
\dfrac{\f}{\f t}\bigg\vert_{t=0}\partial_{k }u_t\circ \varphi_t^{-1} &=-D^2\p(u)\left(\partial_ku,Dg\xi\right)+D\p(u)\left(D(u-g)\partial_k\xi-\partial_kDg\xi\right)\\
&=-D^2\p(u)\left(\partial_ku,Dg\xi\right)+Du\partial_k\xi-D\p(u)\left(Dg\partial_k\xi+\partial_kDg\xi\right).
\numberthis\label{30}
\end{align*}
The second equality in \eqref{30} yields because of $D\p(u)(Du\partial_k\xi)=Du\partial_k\xi$, since $Du\cdot v\in T_u\mathcal{N}$ for a.e. $x\in \Omega$ and all $v\in \R^m$. Analogue to the above computations, we get from equation \eqref{19} with \eqref{22} and \eqref{29}:
\begin{align*}
\dfrac{\f}{\f t}\bigg\vert_{t=0}\partial_{i}\partial_{j}u_t\circ \varphi_t^{-1}&= -D^3\p(u)\left(\partial_iu,\partial_ju,Dg\xi\right)\\
&\quad+D^2\p(u)\left(D(u-g)\partial_i\xi-\partial_iDg\xi,\partial_ju\right)\\
&\quad+D^2\p(u)\left(\partial_iu,D(u-g)\partial_j\xi-\partial_jDg\xi\right)\\
&\quad-  D^2\p(u)\left(\partial_i\partial_ju,Dg\xi\right)+D\p(u)\left(\dfrac{\f}{\f t}\bigg\vert_{t=0}\partial_i\partial_j\U_t\circ \varphi_t^{-1} \right)\numberthis\label{31}
\end{align*}
Due to equation \eqref{16} we have
\begin{align*}
 \partial_i\partial_j\U_t\circ\varphi^{-1}_t&= \partial_k\partial_l(u-g)\left(\partial_i\varphi_t^l\circ\varphi^{-1}_t\right)\left(\partial_j\varphi_t^k\circ\varphi^{-1}_t\right)\\
 &\quad+\left( \partial_k(u-g)\circ \varphi_t\circ\varphi^{-1}_t  \right)\partial_i\partial_j\varphi_t^k\circ\varphi^{-1}_t+\partial_i\partial_jg\circ\varphi^{-1}_t.\numberthis\label{32}
 \end{align*}
Moreover, it holds
\begin{align*}
\dfrac{\f }{\f t}\bigg\vert_{t=0} \partial_i\partial_j\varphi^k_t\circ \varphi^{-1}_t&=\partial_i\partial_j\xi^k.\numberthis\label{33}
\end{align*}
So, we obtain with the equations \eqref{32} and \eqref{33} that
\begin{align*}
\dfrac{\f}{\f t}\bigg\vert_{t=0}\partial_i\partial_j\U_t\circ \varphi_t^{-1}&=\partial_{k}\partial_{l}(u-g)\partial_i\xi^{l}\delta_{jk}+\partial_{k}\partial_{l}(u-g)\partial_j\xi^{k}\delta_{il}+\partial_k(u-g)\partial_i\partial_j\xi^k-\partial_{i}\partial_j\partial_kg\xi^k\\
&=\partial_jD(u-g)\partial_i\xi+\partial_iD(u-g)\partial_j\xi+D(u-g)\partial_i\partial_j\xi-\partial_i\partial_jDg\xi\numberthis\label{34},
\end{align*}
whereby equation \eqref{31} becomes
\begin{align*}
\dfrac{\f}{\f t}\bigg\vert_{t=0}\partial_{i}\partial_{j}u_t\circ \varphi_t^{-1}&= -D^3\p(u)\left(\partial_iu,\partial_ju,Dg\xi\right)\\
&\quad+D^2\p(u)\left(D(u-g)\partial_i\xi-\partial_iDg\xi,\partial_ju\right)\\
&\quad+D^2\p(u)\left(\partial_iu,D(u-g)\partial_j\xi-\partial_jDg\xi\right)\\
&\quad-  D^2\p(u)\left(\partial_i\partial_ju,Dg\xi\right)\\
&\quad+D\p(u)\left(\partial_iD(u-g)\partial_j\xi+\partial_jD(u-g)\partial_i\xi\right)\\
&\quad-D\p(u)\left(Dg\partial_i\partial_j\xi+\partial_i\partial_jDg\xi\right)+Du\partial_i\partial_j\xi\numberthis\label{36}
\end{align*}
where we used $D\pi(u)(Du\partial_i\partial_j\xi)=Du\partial_i\partial_j\xi$ in the last line. Since $Du\cdot v\in T_u\mathcal{N}$ for a.e. $x\in \Omega$ and all $v\in \R^m$, it holds $Du\partial_i\xi=D\p(u)(Du\partial_i\xi)$. Differentiating this with respect to $x^j$, we get
\begin{align*}
\partial_jDu\partial_i\xi+Du\partial_i\partial_j\xi&=D^2\p(u)(Du\partial_i\xi,\partial_ju)+D\p(u)(\partial_jDu\partial_i\xi)\\
&\quad+D\p(u)(Du\partial_i\partial_j\xi).\numberthis\label{38}
\end{align*}
Equation \eqref{38} becomes due to the identity $D\p(u)(Du\partial_i\partial_j\xi)=Du\partial_i\partial_j\xi$,
\begin{align*}
\partial_jDu\partial_i\xi&=D^2\p(u)(Du\partial_i\xi,\partial_ju)+D\p(u)(\partial_jDu\partial_i\xi)\numberthis\label{39}.
\end{align*}
Hence, the equation \eqref{36} reduces to
\begin{align*}
\dfrac{\f}{\f t}\bigg\vert_{t=0}\left(\partial_{i}\partial_{j}u_t\circ \varphi_t^{-1} \right)&= -D^3\p(u)\left(\partial_iu,\partial_ju,Dg\xi\right)\\
&\quad-D^2\p(u)\left(Dg\partial_i\xi+\partial_iDg\xi,\partial_ju\right)-D^2\p(u)\left(\partial_iu,Dg\partial_j\xi+\partial_jDg\xi\right)\\
&\quad-  D^2\p(u)\left(\partial_i\partial_ju,Dg\xi\right)\\
&\quad-D\p(u)\left(\partial_iDg\partial_j\xi+\partial_jDg\partial_i\xi\right)-D\p(u)\left(Dg\partial_i\partial_j\xi+\partial_i\partial_jDg\xi\right)\\
&\quad+\partial_jDu\partial_i\xi+\partial_iDu\partial_j\xi+Du\partial_i\partial_j\xi.\numberthis\label{40}
\end{align*}

\textit{Step 4}. Furthermore, we have because of \eqref{24},
\begin{align*}
\dfrac{\f}{\f t}\bigg\vert_{t=0}\gamma^{ij}\circ \varphi_t^{-1}=-\partial_l\gamma^{ij}\xi^l\text{ and }\dfrac{\f}{\f t}\bigg\vert_{t=0}\Gamma^k_{ij}\circ \varphi_t^{-1} =-\partial_l\Gamma_{ij}^k\xi^l\numberthis\label{41}.
\end{align*}
Putting \eqref{30}, \eqref{40}, and \eqref{41} into \eqref{20} yields the following equation,
\begin{align*}
\dfrac{\f}{\f t}\bigg\vert_{t=0}\Delta_{\gamma}u_t\circ\varphi_t^{-1}&=-\partial_l\gamma^{ij}\xi^l\cdot\left(\partial_i\partial_ju-\Gamma_{ij}^k\partial_ku\right)+\gamma^{ij}\cdot\partial_l\Gamma_{ij}^k\xi^l\partial_ku\\
&\quad+\gamma^{ij}\cdot\left(\partial_jDu\partial_i\xi+\partial_iDu\partial_j\xi+Du\partial_i\partial_j\xi\right)\\
&\quad+\gamma^{ij}\cdot \Gamma^k_{ij}\cdot D^2\p(u)\left(\partial_ku,Dg\xi\right)-Du\gamma^{ij}\cdot \Gamma^k_{ij}\partial_k\xi\\
&\quad+ \gamma^{ij}\cdot \Gamma^k_{ij}\cdot D\p(u)\left(Dg\partial_k\xi+\partial_kDg\xi\right)\\
&\quad-\gamma^{ij}\cdot D^3\p(u)\left(\partial_iu,\partial_ju,Dg\xi\right)\\
&\quad-\gamma^{ij}\cdot D^2\p(u)\left(Dg\partial_i\xi+\partial_iDg\xi,\partial_ju\right)\\
&\quad-\gamma^{ij}\cdot D^2\p(u)\left(\partial_iu,Dg\partial_j\xi+\partial_jDg(x)\xi\right)\\
&\quad-\gamma^{ij}\cdot D^2\p(u)\left(\partial_i\partial_ju,Dg\xi\right)\\
&\quad-\gamma^{ij}\cdot D\p(u)\left(\partial_iDg\partial_j\xi+\partial_jDg\partial_i\xi\right)\\
&\quad-\gamma^{ij}\cdot D\p(u)\left(Dg\partial_i\partial_j\xi+\partial_i\partial_jDg\xi\right).\numberthis\label{43}
\end{align*} 
 
\textit{Step 5}. Now, we continue by determining $\dfrac{\f}{\f t}\bigg\vert_{t=0}\left( \det{D\varphi_t^{-1}} \right)$. For an $(m\times m)$-matrix $(z_{\alpha\beta})_{\alpha,\beta=1}^m$ let be $\text{ad}_{jk}\left((z_{\alpha\beta})_{\alpha,\beta=1}^m\right):=(-1)^{j+k}\det \left((z_{\alpha\beta})_{\alpha,\beta=1}^m\right)_{\alpha\neq j,\beta\neq k}$ the $(m-1)\times (m-1)$-minors for all $1\leq j,k\leq m$. With the Laplacian expansion theorem we deduce $\dfrac{\partial}{\partial z_{jk}}\det((z_{\alpha\beta})_{\alpha,\beta=1}^m)=\text{ad}_{jk}((z_{\alpha\beta})_{\alpha,\beta=1}^m)$. Due to $\text{ad}_{jk}(\text{id})=\delta_{jk}$ we obtain using the chain rule and equation \eqref{24},
\begin{align*}
\dfrac{\f }{\f t}\bigg\vert_{t=0}\text{det }(D\varphi_t^{-1})&=-\div{ \xi}.\numberthis\label{46}
\end{align*}
\textit{Step 6}. Finally, we get by using chain rule and equation \eqref{24} once again,
\begin{align*}
\dfrac{\f}{\f t}\bigg\vert_{t=0}\sqrt{\gamma\circ\varphi_t^{-1}}=-\partial_k\left(\sqrt{\gamma}\right)\xi^k.\numberthis\label{47}
\end{align*}
Now, we put \eqref{43}, \eqref{46} and \eqref{47} into \eqref{10}, and summarize suitably. Then, we obtain
\begin{align*}
\dfrac{1}{\sqrt{\gamma}}\dfrac{\f}{\f t}\bigg\vert_{t=0}f(t,x)&=2 \Delta_{\gamma}u\cdot \partial_l\left(\gamma^{ij}\Gamma_{ij}^k\right)\xi^l\partial_ku-2 \Delta_{\gamma}u \cdot\partial_l\gamma^{ij}\xi^l\partial_i\partial_ju\\
&\quad-\vert \Delta_{\gamma}u\vert^2  \div{\xi}-\vert \Delta_{\gamma}u\vert^2 \dfrac{1}{\sqrt{\gamma}}\partial_k\left(\sqrt{\gamma}\right)\xi^k-2 \Delta_{\gamma}u\cdot\gamma^{ij} \Gamma^k_{ij} Du\partial_k\xi\\
&\quad+2 \Delta_{\gamma}u\cdot\gamma^{ij}\left(\partial_jDu\partial_i\xi+\partial_iDu\partial_j\xi+Du\partial_i\partial_j\xi\right)\\
&\quad+2 \Delta_{\gamma}u\cdot\gamma^{ij}\Gamma^k_{ij} D^2\p(u)\left(\partial_ku,Dg\xi\right)\\
&\quad+ 2 \Delta_{\gamma}u\cdot\gamma^{ij}\Gamma^k_{ij}D\p(u)\left(Dg\partial_k\xi+\partial_kDg\xi\right)\\
&\quad-2 \Delta_{\gamma}u\cdot\gamma^{ij}D^3\p(u)\left(\partial_iu,\partial_ju,Dg\xi\right)\\
&\quad-2 \Delta_{\gamma}u\cdot\gamma^{ij} D^2\p(u)\left(Dg\partial_i\xi+\partial_iDg\xi,\partial_ju\right)\\
&\quad-2 \Delta_{\gamma}u\cdot\gamma^{ij} D^2\p(u)\left(\partial_iu,Dg\partial_j\xi+\partial_jDg\xi\right)\\
&\quad-2 \Delta_{\gamma}u\cdot \gamma^{ij} D^2\p(u)\left(\partial_i\partial_ju,Dg\xi\right)\\
&\quad-2 \Delta_{\gamma}u\cdot\gamma^{ij} D\p(u)\left(\partial_iDg\partial_j\xi+\partial_jDg\partial_i\xi\right)\\
&\quad-2 \Delta_{\gamma}u\cdot \gamma^{ij} D\p(u)\left(Dg\partial_i\partial_j\xi+\partial_i\partial_jDg\xi\right).\numberthis\label{48}
\end{align*} 
Notice that $\div_{\gamma}{\xi}=\div{\xi}+\dfrac{1}{\sqrt{\gamma}}\partial_k\left(\sqrt{\gamma}\right)\xi^k=\div{\xi}+\Gamma^l_{kl}\xi^k$. Furthermore, we can rewrite $2 \Delta_{\gamma}u\cdot\gamma^{ij}\left(\partial_jDu\partial_i\xi+\partial_iDu\partial_j\xi\right)$ into $4\Delta_{\gamma}u\cdot\gamma^{ij}\partial_iDu\partial_j\xi$ and with the Riemannian Gradient $\grad_{\gamma}f=\gamma^{ij}\partial_if\partial_j$ into $4\Delta_{\gamma}u\cdot D^2u\grad_{\gamma}{\xi}$. Moreover, we have $ \gamma^{ij}Du\partial_i\partial_j\xi-\gamma^{ij}\Gamma^k_{ij}Du\partial_k\xi=Du\Delta_{\gamma}\xi$. With this abbreviations and putting \eqref{48} into \eqref{9b} we obtain the \textit{differential equation for variationally biharmonic maps} after reformulations,
\begin{align*}
&\int_{\Omega}\left( 4\Delta_{\gamma}u\cdot D^2u\grad_{\gamma}{\xi}+2\Delta_{\gamma}u\cdot Du\Delta_{\gamma}\xi-\vert \Delta_{\gamma}u\vert^2 \div_{\gamma}{\xi}\right)\sqrt{\gamma}\f\mathcal{L}^m\\
&=\int_{\Omega} 2 \Delta_{\gamma}u\cdot\gamma^{ij} D^3\p(u)\left(\partial_iu,\partial_ju,Dg\xi\right)\sqrt{\gamma}\f\mathcal{L}^m\\
&\quad+\int_{\Omega} 2 \Delta_{\gamma}u\cdot\gamma^{ij} D^2\p(u)\left(Dg\partial_i\xi+\partial_iDg\xi,\partial_ju\right)\sqrt{\gamma}\f\mathcal{L}^m(x)\\
&\quad+\int_{\Omega} 2 \Delta_{\gamma}u\cdot\gamma^{ij} D^2\p(u)\left(\partial_iu,Dg\partial_j\xi+\partial_jDg\xi\right)\sqrt{\gamma}\f\mathcal{L}^m\\
&\quad+\int_{\Omega} 2 \Delta_{\gamma}u\cdot\gamma^{ij} D^2\p(u)\left(\partial_i\partial_ju,Dg\xi\right)\sqrt{\gamma}\f\mathcal{L}^m\\
&\quad+\int_{\Omega} 2 \Delta_{\gamma}u\cdot\gamma^{ij} D\p(u)\left(\partial_iDg\partial_j\xi+\partial_jDg\partial_i\xi\right)\sqrt{\gamma}\f\mathcal{L}^m\\
&\quad+\int_{\Omega} 2 \Delta_{\gamma}u\cdot\gamma^{ij} D\p(u)\left(Dg\partial_i\partial_j\xi+\partial_i\partial_jDg\xi\right)\sqrt{\gamma}\f\mathcal{L}^m\\
&\quad-\int_{\Omega} 2 \Delta_{\gamma}u\cdot\gamma^{ij} \Gamma^k_{ij} D^2\p(u)\left(\partial_ku,Dg\xi\right)\sqrt{\gamma}\f\mathcal{L}^m\\
&\quad-\int_{\Omega} 2 \Delta_{\gamma}u\cdot\gamma^{ij} \Gamma^k_{ij} D\p(u)\left(Dg\partial_k\xi+\partial_kDg\xi\right)\sqrt{\gamma}\f\mathcal{L}^m\\
&\quad-\int_{\Omega} 2 \Delta_{\gamma}u\cdot\partial_l\left(\gamma^{ij}\Gamma_{ij}^k\right)\xi^l\partial_ku\sqrt{\gamma}\f\mathcal{L}^m+\int_{\Omega} 2 \Delta_{\gamma}u\cdot \partial_l\gamma^{ij}\xi^l\partial_i\partial_ju\sqrt{\gamma}\f\mathcal{L}^m\numberthis\label{49}.
\end{align*}
It is straightforward to see that equation \eqref{49} is equivalent to \eqref{7}. This concludes the proof.
\end{proof} 
Notice that equation \eqref{7} takes the form of the equation in Lemma 1 in \cite{1} for Euclidean metric and constant boundary values, since the right-hand side is identical to zero in this case.

\section{Derivation of a boundary monotonicity inequality}\label{sec3}

Before we start with the derivation, we want to mention the following
\begin{remark}
In our estimates we take care to produce 'good-natured' error terms (integrals). We say that an error term is 'good-natured' if the dimension of integration region minus number of derivatives on $u$ is greater than $\vert x\vert$-powers in the denominator. If an error term is good-natured then it vanishes for small radii.
\end{remark}

\begin{proof}[Proof of \Cref{thm}] We derive the boundary monotonicity formula \eqref{57} in 8 steps. All constants appearing in the proof may depend on $m,\mathcal{N},G,H$. Further dependecies will be indicated in parentheses, e.g. $C_1(\Vert Dg\Vert_{L^{\infty}})$.\\
\textit{Step 1.} We set $\Omega=B^+$ in \eqref{49} from the proof of \Cref{lemma1}. Now, we form the right-hand side of \eqref{49} so that no second order derivatives of $\xi$ appear on the right-hand side of \eqref{49}. Moreover, we split $\displaystyle\int_{B^+}2\Delta_{\gamma}u\cdot Du \Delta_{\gamma}\xi\f\mu_{\gamma}$ into $$\int_{B^+}2\Delta_{\gamma}u\cdot Du \gamma^{ij}\partial_i\partial_j\xi\f\mu_{\gamma}-\int_{B^+}2\Delta_{\gamma}u\cdot Du\gamma^{ij}\Gamma^k_{ij}\partial_k\xi\f\mu_{\gamma}$$ and $\displaystyle\int_{B^+}\vert \Delta_{\gamma}u\vert^2 \div_{\gamma}{\xi}\f\mu_{\gamma}$ into $\displaystyle\int_{B^+}\vert \Delta_{\gamma}u\vert^2 \div{\xi}\f\mu_{\gamma}+\int_{B^+}\vert \Delta_{\gamma}u\vert^2 \Gamma^l_{kl}\xi^k\f\mu_{\gamma}$. Then, we bring the second summand respectively on the right-hand side of \eqref{49}. So, we get the following equation,
\begin{align*}
&\int_{B^+}\left( 4\Delta_{\gamma}u\cdot D^2u \grad_{\gamma}{\xi}+2\Delta_{\gamma}u\cdot Du \gamma^{ij}\partial_i\partial_j\xi-\vert \Delta_{\gamma}u\vert^2 \div{\xi}\right)\f\mu_{\gamma}\\
&\quad-\int_{B^+}\left( 2 \Delta_{\gamma}u\cdot D\p(u)\left(Dg\gamma^{ij}\partial_i\partial_j\xi\right)\right)\f\mu_{\gamma}\\
&=\int_{B^+} 2 \Delta_{\gamma}u\cdot\gamma^{ij} D^3\p(u)\left(\partial_iu,\partial_ju,Dg\xi\right)\f\mu_{\gamma}+\int_{B^+} 2 \Delta_{\gamma}u\cdot D^2\p(u)\left(\Delta_{\gamma}u,Dg\xi\right)\f\mu_{\gamma}\\
&\quad+\int_{B^+} 4 \Delta_{\gamma}u\cdot\gamma^{ij} D^2\p(u)\left(\partial_iu,Dg\partial_j\xi+\partial_jDg\xi\right)\f\mu_{\gamma}+\int_{B^+}\vert \Delta_{\gamma}u\vert^2 \Gamma^l_{kl}\xi^k\f\mu_{\gamma}\\
&\quad+\int_{B^+} 4 \Delta_{\gamma}u\cdot\gamma^{ij} D\p(u)\left(\partial_iDg\partial_j\xi\right)\f\mu_{\gamma}+\int_{B^+}2\Delta_{\gamma}u\cdot Du\gamma^{ij}\Gamma^k_{ij}\partial_k\xi\f\mu_{\gamma}\\
&\quad+\int_{B^+} 2 \Delta_{\gamma}u\cdot D\p(u)\left(\Delta_{\gamma}(Dg)\xi\right)\f\mu_{\gamma}-\int_{B^+} 2 \Delta_{\gamma}u\cdot\gamma^{ij} \Gamma^k_{ij} D\p(u)\left(Dg\partial_k\xi\right)\f\mu_{\gamma}\\
&\quad+\int_{B^+} 2 \Delta_{\gamma}u\cdot \partial_l\gamma^{ij}(x)\xi^l\partial_i\partial_ju\f\mu_{\gamma}-\int_{B^+} 2 \Delta_{\gamma}u\cdot\partial_l\left(\gamma^{ij}\Gamma_{ij}^k\right)\xi^l\partial_ku\f\mu_{\gamma}\numberthis\label{58}.
\end{align*}
\textit{Step 2}. Next, we estimate the left-hand side of \eqref{58} by the right-hand side of \eqref{58} and abbreviate the 'left-hand side of' by $LHS$. We obtain,
\begin{align*}
\vert LHS\eqref{58}\vert&\leq 2C_1\int_{B^+}\vert \Delta_{\gamma}u\vert\vert Du\vert^2    \vert\xi\vert  \f \mu_{\gamma}+  C_2\int_{B^+}\vert \Delta_{\gamma}u\vert^2    \vert\xi\vert  \f \mu_{\gamma}+ 4C_3\int_{B^+}\vert \Delta_{\gamma}u\vert \vert Du\vert    \vert D\xi\vert  \f \mu_{\gamma} \\
&\quad+4C_4\int_{B^+}\vert \Delta_{\gamma}u\vert \vert Du\vert    \vert \xi\vert  \f \mu_{\gamma}+C_5\int_{B^+}\vert \Delta_{\gamma}u\vert^2    \vert\xi\vert  \f \mu_{\gamma}+4C_6\int_{B^+}\vert \Delta_{\gamma}u\vert    \vert D\xi\vert  \f \mu_{\gamma}\\
&\quad+2C_7\int_{B^+}\vert \Delta_{\gamma}u\vert \vert Du\vert    \vert D\xi\vert  \f \mu_{\gamma}
+2C_8\int_{B^+}\vert \Delta_{\gamma}u\vert  \vert \xi\vert  \f \mu_{\gamma}+2C_9\int_{B^+}\vert \Delta_{\gamma}u\vert    \vert D\xi\vert  \f \mu_{\gamma}\\
&\quad+2C_{10}\int_{B^+}\vert \Delta_{\gamma}u\vert\vert D^2u\vert    \vert \xi\vert  \f \mu_{\gamma}+2C_{11}\int_{B^+}\vert \Delta_{\gamma}u\vert\vert Du\vert    \vert \xi\vert  \f \mu_{\gamma} \numberthis\label{59},
 \end{align*}
where $C_1=C_1(\Vert Dg\Vert_{\infty})$, $C_2=C_2(\Vert Dg\Vert_{\infty})$, $C_3=C_3(\Vert Dg\Vert_{\infty})$, $C_4=C_4(\Vert Dg\Vert_{C^1})$, $C_6=C_6(\Vert Dg\Vert_{C^1})$, $C_8=C_8(\Vert Dg\Vert_{C^2})$ and $C_9=C_9(\Vert Dg\Vert_{\infty})$.

For all $\nu\in \N$ we choose a function $\psi_{\nu}\in C^{\infty}(\R_+,[0,1])$ with $\psi_{\nu}\equiv1$ on $[0,1-\frac{1}{\nu}]$, $\psi_{\nu}\equiv0$ on $[1,\infty)$, $\psi_{\nu}'\leq0$ and $\displaystyle\int_{\R_+}\vert \psi_{\nu}'\vert=1$. Thereby, we define for $0<\tau<1$ and $a\in \R^{m-1}\times \left\lbrace 0\right\rbrace $ the cut-off function $\xi(x):=\si(x)\cdot (x-a)=\psi_{\nu}\left(\dfrac{\vert x-a\vert}{\tau}\right)\cdot (x-a)$. We assume without loss of generality that $a=0$. Thus, we have $\vert \xi\vert\leq\vert x\vert \si $ and $\vert D\xi\vert \leq \dfrac{\vert x\vert}{\tau}\vert \si'(x)\vert+\si(x)$ where $\si'(x):=\psi_{\nu}'\left(\dfrac{\vert x\vert}{\tau}\right)$. Therefore, we get for the right-hand side of \eqref{59},
\begin{align*}
&RHS\eqref{59}\\
&\leq  2C_1(\Vert Dg\Vert_{\infty})\int_{B^+}\vert \Delta_{\gamma}u\vert\vert Du\vert^2 \vert x\vert   \si  \f \mu_{\gamma}+C_{12}(\Vert Dg\Vert_{\infty})\int_{B^+}\vert \Delta_{\gamma}u\vert^2   \vert x\vert   \si  \f \mu_{\gamma}\\
&\quad+ 4C_{13}(\Vert Dg\Vert_{\infty})\int_{B^+}\vert \Delta_{\gamma}u\vert \vert Du\vert    \si \f \mu_{\gamma}+ \dfrac{4}{\tau} C_{13}(\Vert Dg\Vert_{\infty})\int_{B^+}\vert \Delta_{\gamma}u\vert \vert Du\vert    \vert x\vert \vert\si'\vert  \f \mu_{\gamma}\\
&\quad+4C_{14}(\Vert Dg\Vert_{C^1})\int_{B^+}\vert \Delta_{\gamma}u\vert \vert Du\vert    \vert x\vert   \si  \f \mu_{\gamma}+2C_{15}(\Vert Dg\Vert_{C^1})\int_{B^+}\vert \Delta_{\gamma}u\vert    \si  \f \mu_{\gamma}\\
&\quad+\dfrac{2}{\tau} C_{15}(\Vert Dg\Vert_{C^1})\int_{B^+}\vert \Delta_{\gamma}u\vert   \vert x\vert \vert \si'\vert  \f \mu_{\gamma}+2C_8(\Vert Dg\Vert_{C^2})\int_{B^+}\vert \Delta_{\gamma}u\vert  \vert x\vert   \si  \f \mu_{\gamma}\\
&\quad+2C_{10}\int_{B^+}\vert \Delta_{\gamma}u\vert\vert D^2u\vert  \vert x\vert   \si  \f \mu_{\gamma}\\
&=\int_{B^+}\left(2C_1\vert \Delta_{\gamma}u\vert\vert Du\vert^2 +4C_{14}\vert \Delta_{\gamma}u\vert \vert Du\vert+2C_8\vert \Delta_{\gamma}u\vert +2C_{10}\vert \Delta_{\gamma}u\vert\vert D^2u\vert+C_{12}\vert \Delta_{\gamma}u\vert^2\right)\vert x\vert   \si  \f \mu_{\gamma}\\
&\quad+\int_{B^+}\left(4C_{13}\vert \Delta_{\gamma}u\vert \vert Du\vert +2C_{15}\vert \Delta_{\gamma}u\vert\right) \si  \f \mu_{\gamma}\\
&\quad+\dfrac{1}{\tau} \int_{B^+}\left(4C_{13}\vert \Delta_{\gamma}u\vert \vert Du\vert +2C_{15} \vert \Delta_{\gamma}u\vert\right)  \vert x\vert \vert\si'\vert  \f \mu_{\gamma}\numberthis\label{60}\\
&=:I+II+III,
\end{align*}
where $C_{12}=C_{12}(\Vert Dg\Vert_{\infty}):=C_2(\Vert Dg\Vert_{\infty})+C_5 $, $C_{13}=C_{13}(\Vert Dg\Vert_{\infty}):=C_3(\Vert Dg\Vert_{\infty})+\dfrac{1}{2} C_7$, $C_{14}=C_{14}(\Vert Dg\Vert_{C^1}):=C_4(\Vert Dg\Vert_{C^1})+\frac{1}{2}C_{11}$ and $C_{15}=C_{15}(\Vert Dg\Vert_{C^1}):=2C_6(\Vert Dg\Vert_{C^1})+C_9(\Vert Dg\Vert_{\infty})$. Now, we estimate $I$, $II$ and $III$ as follows. With the hepl of Young's inequality, we estimate
\begin{align*}
I&\leq C_{16}\int_{B^+}\vert \Delta_{\gamma}u\vert^2\vert x\vert   \si  \f \mu_{\gamma}+C_{17}\int_{B^+}\vert Du\vert^4 \vert x\vert   \si  \f \mu_{\gamma} + C_{10}\int_{B^+}\vert D^2u\vert^2 \vert x\vert   \si  \f \mu_{\gamma}\\
&\quad+C_{18}\int_{B^+}\vert x\vert   \si  \f \mu_{\gamma} \numberthis\label{61}
\end{align*}
where $C_{16}=C_{16}(\Vert Dg\Vert_{C^2}):=C_1+C_8+C_{10}+C_{12}+2C_{14}$, $C_{17}=C_{17}(\Vert Dg\Vert_{C^1}):=C_1+C_{14}$ and $C_{18}=C_{18}(\Vert Dg\Vert_{C^2}):=C_8+C_{14}$. Further, we obtain
\begin{align*}
&II\\
&\leq C_{19}\int_{B^+}\vert \Delta_{\gamma}u\vert^2\vert x\vert   \si  \f \mu_{\gamma}+C_{13}\int_{B^+}\vert Du\vert^4\vert x\vert   \si  \f \mu_{\gamma}+\int_{B^+}\left(\dfrac{C_{13}}{\vert x\vert}+\dfrac{C_{15}}{\vert x\vert^3}\right)   \si  \f \mu_{\gamma}\numberthis\label{62}
\end{align*}
with $C_{19}=C_{19}(\Vert Dg\Vert_{C^1}):=2C_{13}+C_{15}$. Moreover, we get due to $\vert \Delta_{\gamma}u\vert \leq G\vert D^2u\vert+C_5\vert Du\vert$ and applying Young's inequality,
\begin{align*}
&III\\
&\leq \dfrac{1}{\tau} \int_{B^+}\left(C_{20}\vert D^2u\vert^2\vert x\vert+2C_{13}G\dfrac{\vert Du\vert^2}{\vert x\vert}+C_{21}\vert Du\vert^2+\dfrac{C_{15}G}{\vert x\vert}+C_{15}C_5 \right)  \vert x\vert \vert\si'\vert  \f \mu_{\gamma},\numberthis\label{63}
\end{align*}
where $C_{20}=C_{20}(\Vert Dg\Vert_{C^1}):=(2C_{13}+C_{15})G$ and $C_{21}=C_{21}(\Vert Dg\Vert_{C^1}):=(4C_{13}+C_{15})C_5$. Together with \eqref{61}, \eqref{62} and \eqref{63} we obtain, since $\vert x\vert \leq 1$ on the domain of integration
\begin{align*}
&\vert LHS\eqref{58}\vert\leq I+II+III\\
&\leq  \int_{B^+}\left(C_{22}\vert \Delta_{\gamma}u\vert^2+C_{23}\vert Du\vert^4+C_{10}\vert D^2u\vert^2\right)\vert x\vert   \si  \f \mu_{\gamma}+\int_{B^+}\left(\dfrac{C_{13}+C_{15}+C_{18}}{\vert x\vert^3}\right)   \si  \f \mu_{\gamma}\\
&\quad+\dfrac{1}{\tau} \int_{B^+}\left(C_{20}\vert D^2u\vert^2\vert x\vert+(2C_{13}G+C_{21})\dfrac{\vert Du\vert^2}{\vert x\vert}+\dfrac{C_{15}(G+C_5)}{\vert x\vert} \right)  \vert x\vert \vert\si'\vert  \f \mu_{\gamma}\numberthis\label{64}\\
&=:IV+V+VI
\end{align*}
with $C_{22}=C_{22}(\Vert Dg\Vert_{C^2}):=C_{16}+C_{19}$ and $C_{23}=C_{23}(\Vert Dg\Vert_{C^1}):=C_{13}+C_{17}$. Now, we multiply the inequality \eqref{64} by $\e^{\chi\tau}\tau^{3-m}$ and integrate over $[\rho,r]$,
\begin{align*}
\int_{\rho}^r\e^{\chi\tau}\tau^{3-m}\vert LHS(58)\vert\f \tau\leq \int_{\rho}^r\e^{\chi\tau}\tau^{3-m}(IV+V+VI)\f \tau\numberthis\label{65},
\end{align*}
where $0<\rho<r<R$. $\si(x)$ converge to the characteristic function of $B^+_{\tau}$ as $\nu\rightarrow\infty$. Thus, applying the dominated convergence theorem and estimating $\vert x\vert<\tau$, we obtain
\begin{align*}
\lim_{\nu\rightarrow\infty}\int_{\rho}^r\e^{\chi\tau}\tau^{3-m}(IV+V)\f \tau&\leq \int_{\rho}^r\e^{\chi\tau}\tau^{4-m}\int_{B^+_{\tau}}\left(C_{22}\vert \Delta_{\gamma}u\vert^2+C_{23}\vert Du\vert^4+C_{10}\vert D^2u\vert^2\right) \f\mu_{\gamma}\f\tau\\
&\quad+\tilde{C}_{24}\int_{\rho}^r\e^{\chi\tau}\tau^{3-m}\int_{B^+_{\tau}}\dfrac{1}{\vert x\vert^3}  \f\mu_{\gamma}\f\tau\numberthis\label{66},
\end{align*}
for a.e. $0<\rho<r< R$ where $\tilde{C}_{24}:=\tilde{C}_{24}(\Vert Dg\Vert_{C^2})$. The square roots of the eigenvalues $\lambda_l$ of $(\gamma_{ij})$ lie in $(G^{-1/2},G^{1/2})$. Hence, $\sqrt{\gamma}=\sqrt{\det{\gamma_{ij}(x)}}=\displaystyle\prod_{l=1}^m\lambda_l^{1/2}\leq G^{m/2}$. Thus, it follows for the last integral in \eqref{66},
\begin{align*}
\tilde{C}_{24}\int_{\rho}^r\e^{\chi\tau}\tau^{3-m}\int_{B^+_{\tau}}\dfrac{1}{\vert x\vert^3}  \f\mu_{\gamma}\f\tau&\leq C_{24}\cdot (r-\rho) \numberthis\label{68}
\end{align*}
with $C_{24}=C_{24}(\Vert Dg\Vert_{C^2},\Vert u\Vert_{L^{\infty}})$. Furthermore, it holds by Gagliardo-Nirenberg's  interpolation inequality (cf. \cite{21}, page 125 and 126) for $j=1$, $p=4$, $k=r=2$, $\alpha=\dfrac{1}{2}$ and $q,s=\infty$
\begin{align*}
\Vert Dv\Vert_{L^4}\leq \tilde{C}_1 \Vert D^2v\Vert_{L^2}^{1/2}\Vert v\Vert_{L^{\infty}}^{1/2}+\tilde{C}_2\Vert v\Vert_{L^{\infty}}\numberthis\label{70}.
\end{align*}
Rescaling from $B_1^+$ to $B_{\tau}^+$, we obtain the following version of this estimate:
\begin{align*}
\tau^{4-m}\int_{\B}\vert Du\vert^4\f \mu_{\gamma}&\leq C_{25}(B^+_1)\Vert u\Vert^2_{L^{\infty}(B^+_1)}\tau^{4-m}\int_{\B}\vert D^2u\vert^2\f \mu_{\gamma}+C_{26}(B^+_1)\Vert u\Vert^4_{L^{\infty}(B^+_1)}\numberthis\label{72}
\end{align*}
For the right-hand side of \eqref{66} we get with \eqref{68} and \eqref{72} the following estimate,
\begin{align*}
RHS\eqref{66}&\leq C_{22} \int_{\rho}^r\e^{\chi\tau}\tau^{4-m}\int_{B^+_{\tau}}\vert \Delta_{\gamma}u\vert^2\f\mu_{\gamma}\f\tau+C_{27}\int_{\rho}^r\e^{\chi\tau}\tau^{4-m}\int_{B^+_{\tau}}\vert D^2u\vert^2 \f\mu_{\gamma}\f\tau\\
&\quad+C_{28}\cdot  (r-\rho)\numberthis\label{73}
\end{align*}
where $C_{27}=C_{27}(\Vert Dg\Vert_{C^1},B^+_1,\Vert u\Vert_{L^{\infty}(B_1^+)}):=C_{10}+C_{23}\cdot C_{25}\Vert u\Vert^2_{L^{\infty}(B_1^+)}$ and\\ $C_{28}=C_{28}(\Vert Dg\Vert_{C^2},\Vert u\Vert_{L^{\infty}(B_1^+)} ):=C_{24}+C_{23}\cdot C_{26}\Vert u\Vert^4_{L^{\infty}(B_1^+)}\e^{\chi R}$. Thanks to Lemma 2 in the appendix of \cite{1} we obtain for a.e. $0<\rho<r<R$:
\begin{align*}
&\lim_{\nu\rightarrow\infty}\int_{\rho}^r\e^{\chi\tau}\tau^{3-m}VI\f \tau\\
&=\int_{B^+_r\setminus B^+_{\rho}}\e^{\chi \vert x\vert}\left(C_{20}\dfrac{\vert D^2u\vert^2}{\vert x\vert^{m-5}} +(2C_{13}G+C_{21})\dfrac{\vert Du\vert^2}{\vert x\vert^{m-3}}  +\dfrac{C_{15}(G+C_5)}{\vert x\vert^{m-3}}\right) \f \mu_{\gamma}\\
&\leq \int_{B^+_r\setminus B^+_{\rho}}\e^{\chi \vert x\vert}\left(C_{20}\dfrac{\vert D^2u\vert^2}{\vert x\vert^{m-5}} +(2C_{13}G+C_{21})\dfrac{\vert Du\vert^2}{\vert x\vert^{m-3}}\right) \f \mu_{\gamma}+C_{29}\cdot (r-\rho)\numberthis\label{78}
\end{align*}
with $C_{29}=C_{29}(\Vert Dg\Vert_{C^2},\Vert u\Vert_{L^{\infty}}):=G^{m/2}\h^{m-1}(S_1^+)\e^{\chi R}C_{15}(G+C_{5})R^2$. Together with \eqref{73} and \eqref{78} we have,
\begin{align*}
&\lim_{\nu\rightarrow\infty} \int_{\rho}^r\e^{\chi\tau}\tau^{3-m}\vert LHS\eqref{58}\vert \f \tau\leq \lim_{\nu\rightarrow\infty} \int_{\rho}^r\e^{\chi\tau}\tau^{3-m}(IV+V+VI) \f \tau\\
&\leq  C_{22} \int_{\rho}^r\e^{\chi\tau}\tau^{4-m}\int_{B^+_{\tau}}\vert \Delta_{\gamma}u\vert^2\f\mu_{\gamma}\f\tau+C_{27}\int_{\rho}^r\e^{\chi\tau}\tau^{4-m}\int_{B^+_{\tau}}\vert D^2u\vert^2 \f\mu_{\gamma}\f\tau+C_{30}\cdot (r-\rho)\\
&\quad+\int_{B^+_r\setminus B^+_{\rho}}\e^{\chi \vert x\vert}\left(C_{20}\dfrac{\vert D^2u\vert^2}{\vert x\vert^{m-5}} +C_{31}\dfrac{\vert Du\vert^2}{\vert x\vert^{m-3}}\right) \f \mu_{\gamma}\numberthis\label{79}
\end{align*}
where $C_{30}=C_{30}(\Vert Dg\Vert_{C^2},\Vert u\Vert_{L^{\infty}}):=C_{28}+C_{29}$ and $C_{31}=C_{31}(\Vert Dg\Vert_{C^1}):=2C_{13}G+C_{21}$.\\

\textit{Step 3.} Notice that it holds
\begin{align*}
\lim_{\nu\rightarrow\infty}\bigg\vert \int_{\rho}^r\e^{\chi\tau}\tau^{3-m} LHS\eqref{58} \f \tau\bigg\vert \leq\lim_{\nu\rightarrow\infty} \int_{\rho}^r\e^{\chi\tau}\tau^{3-m}\vert LHS\eqref{58}\vert \f \tau\leq RHS\eqref{79}\numberthis\label{80}
\end{align*} 
because of the monotonicity for integrals. We will find an estimate for $\displaystyle\int_{\rho}^r\e^{\chi\tau}\tau^{3-m} LHS\eqref{58} \f \tau.$ For that purpose we rewrite $LHS\eqref{58}$ as follows,
\begin{align*}
LHS\eqref{58}&=\int_{B^+}\left(4\Delta_{\gamma}u\cdot\gamma^{ij}\partial_iDu\partial_j\xi+2\Delta_{\gamma}u\cdot Du\Delta\xi-\vert \Delta_{\gamma}u\vert^2\div{\xi}\right)\f\mu_{\gamma}\\
&\quad+2\int_{B^+}\left(\Delta_{\gamma}u\cdot Du(\gamma^{ij}-\delta_{ij})\partial_i\partial_j\xi -\Delta_{\gamma}u\cdot D\p(u)\left( Dg\cdot\gamma^{ij}\partial_i\partial_j\xi   \right)\right)\f \mu_{\gamma}\numberthis\label{81}\\
&=:VII+VIII
\end{align*}
For the sake of clarity we use $f_i$ for partial derivatives $\partial_if$ and write $f_{ij}$ instead of $\partial_i\partial_jf$. Putting $\xi^k_j(x)=\dfrac{1}{\tau}\si'\dfrac{x^jx^k}{\vert x\vert}+\si\delta_{jk}$, $\div{\xi(x)}=\dfrac{\vert x\vert}{\tau}\si'+m\si$ and $\Delta \xi^k(x)=\dfrac{1}{\tau^2}\si''x^k+\dfrac{m+1}{\tau}\si'\dfrac{x^k}{\vert x\vert}$ where $\si:=\psi_{\nu}\left(\dfrac{\vert x\vert}{\tau}\right)$, $\si':=\psi_{\nu}'\left(\dfrac{\vert x\vert}{\tau}\right)$ and $\si'':=\psi_{\nu}''\left(\dfrac{\vert x\vert}{\tau}\right)$ into $VII$ we obtain
\begin{align*}
VII&=4\int_{B^+}\Delta_{\gamma}u\cdot\gamma^{ik}u_{ik} \si\f\mu_{\gamma}-m\int_{B^+}\vert \Delta_{\gamma}u\vert^2 \si\f\mu_{\gamma}\\
&\quad+\dfrac{4}{\tau}\int_{B^+}\dfrac{\Delta_{\gamma}u\cdot\gamma^{ij}u_{ik}x^jx^k}{\vert x\vert} \si'\f\mu_{\gamma}+\dfrac{2(m+1)}{\tau}\int_{B^+}\dfrac{\Delta_{\gamma}u\cdot u_k x^k}{\vert x\vert} \si'\f\mu_{\gamma}\\
&\quad+\dfrac{2}{\tau^2}\int_{B^+}\Delta_{\gamma}u\cdot u_k x^k \si''\f\mu_{\gamma}-\dfrac{1}{\tau}\int_{B^+}\vert \Delta_{\gamma}u\vert^2\vert x\vert \si'\f\mu_{\gamma}\numberthis\label{82}
\end{align*}
Using the Laplace-Beltrami-operator $\Delta_{\gamma}u=\gamma^{ik}u_{ik}-\gamma^{ik}\Gamma_{ik}^lu_{l}$ we can rewrite the first integral in \eqref{82} as follows,
\begin{align*}
4\int_{B^+}\Delta_{\gamma}u\cdot\gamma^{ik}u_{ik} \si\f\mu_{\gamma}=4\int_{B^+}\vert\Delta_{\gamma}u\vert^2\si\f\mu_{\gamma} +4\int_{B^+}\Delta_{\gamma}u\cdot\gamma^{ik}\Gamma_{ik}^lu_{l} \si\f\mu_{\gamma}\numberthis\label{83}
\end{align*}
Moreover, we rewrite the third integral in \eqref{82} with $\gamma^{ij}=\delta^{ij}+(\gamma^{ij}-\delta^{ij})$ as follows,
\begin{align*}
\dfrac{4}{\tau}\int_{B^+}\dfrac{\Delta_{\gamma}u\cdot\gamma^{ij}u_{ik}x^jx^k}{\vert x\vert} \si'\f\mu_{\gamma}&=\dfrac{4}{\tau}\int_{B^+}\dfrac{\Delta_{\gamma}u\cdot u_{jk}x^jx^k}{\vert x\vert} \si'\f\mu_{\gamma}\\
&\quad+\dfrac{4}{\tau}\int_{B^+}\dfrac{\Delta_{\gamma}u\cdot(\gamma^{ij}-\delta^{ij})u_{ik}x^jx^k}{\vert x\vert} \si'\f\mu_{\gamma}\numberthis\label{84}
\end{align*}
So, we have with \eqref{83} and \eqref{84} for \eqref{82},
\begin{align*}
VII&=(4-m)\int_{B^+}\vert \Delta_{\gamma}u\vert^2 \si\f\mu_{\gamma}-\dfrac{1}{\tau}\int_{B^+}\vert \Delta_{\gamma}u\vert^2\vert x\vert \si'\f\mu_{\gamma}+\dfrac{4}{\tau}\int_{B^+}\dfrac{\Delta_{\gamma}u\cdot u_{jk}x^jx^k}{\vert x\vert} \si'\f\mu_{\gamma}\\
&\quad+\dfrac{2(m+1)}{\tau}\int_{B^+}\dfrac{\Delta_{\gamma}u\cdot u_k x^k}{\vert x\vert} \si'\f\mu_{\gamma}+\dfrac{2}{\tau^2}\int_{B^+}\Delta_{\gamma}u\cdot u_k x^k \si''\f\mu_{\gamma}\\
&\quad+4\int_{B^+}\Delta_{\gamma}u\cdot\gamma^{ik}\Gamma_{ik}^lu_{l} \si\f\mu_{\gamma}+\dfrac{4}{\tau}\int_{B^+}\dfrac{\Delta_{\gamma}u\cdot(\gamma^{ij}-\delta_{ij})u_{ik}x^jx^k}{\vert x\vert} \si'\f\mu_{\gamma}\numberthis\label{85}
\end{align*}
Multiplying $VII$ by  $\e^{\chi\tau}\tau^{3-m}$ and integrating over $[\rho,r]$ yields
\begin{align*}
&\int_{\rho}^r\e^{\chi\tau}\tau^{3-m}VII\f\tau\\
&=\int_{\rho}^r\e^{\chi\tau}(4-m)\tau^{3-m}\int_{B^+}\vert \Delta_{\gamma}u\vert^2 \si\f\mu_{\gamma}\f\tau-\int_{\rho}^r\e^{\chi\tau}\tau^{2-m}\int_{B^+}\vert \Delta_{\gamma}u\vert^2 \vert x\vert\si'\f\mu_{\gamma}\f\tau\\
&\quad+\int_{\rho}^r\e^{\chi\tau}\tau^{2-m}\int_{B^+}\dfrac{4\Delta_{\gamma}u\cdot u_{jk}x^jx^k}{\vert x\vert} \si'\f\mu_{\gamma}\f\tau\\
&\quad+2(m+1)\int_{\rho}^r\e^{\chi\tau}\tau^{2-m}\int_{B^+}\dfrac{\Delta_{\gamma}u\cdot u_k x^k}{\vert x\vert} \si'\f\mu_{\gamma}\f\tau\\
&\quad+2\int_{\rho}^r\e^{\chi\tau}\tau^{1-m}\int_{B^+}\Delta_{\gamma}u\cdot u_k x^k\si''\f\mu_{\gamma}\f\tau+4\int_{\rho}^r\e^{\chi\tau}\tau^{3-m}\int_{B^+}\Delta_{\gamma}u\cdot\gamma^{ik}\Gamma_{ik}^lu_{l} \si\f\mu_{\gamma}\f\tau\\
&\quad+4\int_{\rho}^r\e^{\chi\tau}\tau^{2-m}\int_{B^+}\dfrac{\Delta_{\gamma}u\cdot(\gamma^{ij}-\delta_{ij})u_{ik}x^jx^k}{\vert x\vert} \si'\f\mu_{\gamma}\f\tau\numberthis\label{86}.
\end{align*}
We set $I_{\nu}(\tau):=\displaystyle\tau^{4-m}\int_{\R^m_+}\vert \Delta_{\gamma} u\vert^2\si\f x$. It holds
\begin{align*}
I_{\nu}'(\tau)&=(4-m)\tau^{3-m}\int_{B^+}\vert \Delta_{\gamma} u\vert^2\si\f \mu_{\gamma}-\tau^{2-m}\int_{B^+}\vert \Delta_{\gamma} u\vert^2\vert x\vert\si'\f \mu_{\gamma}\numberthis\label{87}.
\end{align*}
According to Lemma 2 in the appendix of \cite{1} and the dominated convergence theorem, $\displaystyle\int_{\rho}^r\e^{\chi \tau} I_{\nu}'(\tau)\f\tau$ tends for $\nu\rightarrow\infty$  to
\begin{align*}
\int_{\rho}^r\e^{\chi \tau}\left((4-m)\tau^{3-m}\int_{B^+_{\tau}}\vert \Delta_{\gamma} u\vert^2\f x+\tau^{4-m}\int_{\partial B^+_{\tau}}\vert \Delta_{\gamma} u\vert^2\f x\right)\f\tau=\int_{\rho}^r\e^{\chi \tau}I'(\tau)\f\tau\numberthis\label{88}
\end{align*} 
for all $\rho,r$. Since $\frac{\f}{\f\tau}(\si')=-\frac{1}{\tau^2}\si''\vert x\vert$ we get by using Fubini's theorem and applying integration by parts
\begin{align*} 
&2\int_{\rho}^r\e^{\chi\tau}\tau^{1-m}\int_{B^+}\Delta_{\gamma}u\cdot u_k x^k\si''\f\mu_{\gamma}\f\tau\\
&=-2\int_{\rho}^r\e^{\chi\tau}\tau^{3-m}\int_{B^+}\Delta_{\gamma}u\cdot u_k x^k\dfrac{\f}{\f\tau}(\si')\dfrac{1}{\vert x\vert}\f\mu_{\gamma}\f\tau\\
&=-2\int_{B^+}\Delta_{\gamma}u\cdot u_k x^k\left(\int_{\rho}^r\e^{\chi\tau}\tau^{3-m}\dfrac{\f}{\f\tau}(\si')\f\tau\right)\dfrac{1}{\vert x\vert}\f\mu_{\gamma}\\
&=-2\int_{B^+}\Delta_{\gamma}u\cdot u_k x^k\left(\e^{\chi r}r^{3-m}\psi'_{\nu,r}-\e^{\chi \rho}\rho^{3-m}\psi'_{\nu,\rho}\right)\dfrac{1}{\vert x\vert}\f\mu_{\gamma}\\
&\quad+2(3-m)\int_{\rho}^r\e^{\chi\tau}\tau^{2-m}\int_{B^+}\Delta_{\gamma}u\cdot u_k x^k\si'\dfrac{1}{\vert x\vert}\f\mu_{\gamma}\f\tau\\
&\quad+2\chi\int_{\rho}^r\e^{\chi\tau}\tau^{3-m}\int_{B^+}\Delta_{\gamma}u\cdot u_k x^k \si'\dfrac{1}{\vert x\vert}\f\mu_{\gamma}\f\tau\numberthis\label{89}
\end{align*}
for a.e. $0<\rho<r<R$. Furthermore, it holds for the last two integrals in \eqref{86},
\begin{align*}
&4\int_{\rho}^r\e^{\chi\tau}\tau^{3-m}\int_{B^+}\Delta_{\gamma}u\cdot\gamma^{ik}\Gamma_{ik}^lu_{l} \si\f\mu_{\gamma}\f\tau\\
&\quad+4\int_{\rho}^r\e^{\chi\tau}\tau^{2-m}\int_{B^+}\dfrac{\Delta_{\gamma}u\cdot(\gamma^{ij}-\delta_{ij})u_{ik}x^jx^k}{\vert x\vert} \si'\f\mu_{\gamma}\f\tau\\
&\leq 4C_7\int_{\rho}^r\e^{\chi\tau}\tau^{3-m}\int_{B^+}\vert \Delta_{\gamma}u\vert \vert Du\vert \si\f\mu_{\gamma}\f\tau\\
&\quad+4H\int_{\rho}^r\e^{\chi\tau}\tau^{2-m}\int_{B^+}\vert \Delta_{\gamma}u\vert \vert D^2u\vert \vert x\vert^2 \vert \si'\vert\f\mu_{\gamma}\f\tau\numberthis\label{90}.
\end{align*}
Here, we have used the inequality 
\begin{align*}
\vert \gamma^{ij}(x)-\delta^{ij}\vert \leq \vert D\gamma^{ij}\vert \vert x\vert\leq H\vert x\vert\numberthis\label{91}
\end{align*}
which follows from the assumption that $\gamma_{ij}(0)=\delta_{ij}$. We obtain with \eqref{87}, \eqref{89} and \eqref{90} the following estimate for the right-hand side of \eqref{86},
\begin{align*}
RHS\eqref{86}&\leq \int_{\rho}^r\e^{\chi\tau}I_{\nu}'(\tau)\f\tau+\int_{\rho}^r\e^{\chi\tau}\tau^{2-m}\int_{B^+}\left(\dfrac{4\Delta_{\gamma}u\cdot u_{jk}x^jx^k}{\vert x\vert}+\dfrac{8\Delta_{\gamma}u\cdot u_k x^k}{\vert x\vert}\right) \si'\f\mu_{\gamma}\f\tau\\
&\quad-2\int_{B^+}\Delta_{\gamma}u\cdot u_k x^k\left(\e^{\chi r}r^{3-m}\psi'_{\nu,r}-\e^{\chi \rho}\rho^{3-m}\psi'_{\nu,\rho}\right)\dfrac{1}{\vert x\vert}\f\mu_{\gamma}\\
&\quad+2\chi\int_{\rho}^r\e^{\chi\tau}\tau^{3-m}\int_{B^+}\Delta_{\gamma}u\cdot u_k x^k \si'\dfrac{1}{\vert x\vert}\f\mu_{\gamma}\f\tau\\
&\quad+4C_7\int_{\rho}^r\e^{\chi\tau}\tau^{3-m}\int_{B^+}\vert \Delta_{\gamma}u\vert \vert Du\vert \si\f\mu_{\gamma}\f\tau\\
&\quad+4H\int_{\rho}^r\e^{\chi\tau}\tau^{2-m}\int_{B^+}\vert \Delta_{\gamma}u\vert \vert D^2u\vert \vert x\vert^2 \vert \si'\vert\f\mu_{\gamma}\f\tau\numberthis\label{92}.
\end{align*}
Thanks to Lemma 2 in the appendix of \cite{1}, the dominated convergence theorem and Lebesgue's differentiation theorem we obtain together with \eqref{88},
\begin{align*}
&\lim_{\nu\rightarrow\infty}RHS\eqref{92}\\
&=\int_{\rho}^r\e^{\chi\tau}I'(\tau)\f\tau-2\chi\int_{B^+_r\setminus B_{\rho}^+}\e^{\chi\vert x\vert}\dfrac{\Delta_{\gamma}u\cdot u_k x^k}{\vert x\vert^{m-3}} \f\mu_{\gamma}\\
&\quad+\int_{S^+_r\setminus S_{\rho}^+}\e^{\chi\vert x\vert}\dfrac{2\Delta_{\gamma}u\cdot u_k x^k}{\vert x\vert^{m-3}} \sqrt{\gamma}\f \h^{m-1}-\int_{B^+_r\setminus B_{\rho}^+}\e^{\chi\vert x\vert}\left(\dfrac{4\Delta_{\gamma}u\cdot u_{jk}x^jx^k}{\vert x\vert^{m-2}}+\dfrac{8\Delta_{\gamma}u\cdot u_k x^k}{\vert x\vert^{m-2}}\right) \f\mu_{\gamma}\\
&\quad+4C_7\int_{\rho}^r\e^{\chi\tau}\tau^{3-m}\int_{B^+_{\tau}}\vert \Delta_{\gamma}u\vert \vert Du\vert \f\mu_{\gamma}\f\tau+4H\int_{B^+_r\setminus B_{\rho}^+}\e^{\chi\vert x\vert}\dfrac{\vert \Delta_{\gamma}u\vert \vert D^2u\vert}{ \vert x\vert^{m-5}} \f\mu_{\gamma}\numberthis\label{93}.
\end{align*}
for a.e. $0<\rho<r<R$. With $2\vert \Delta_{\gamma}u\vert \vert Du\vert\leq \tau \vert \Delta_{\gamma}u\vert^2+\dfrac{1}{\tau} \vert Du\vert^2$ we estimate the second to last integral in \eqref{93} as follows,
\begin{align*}
&4C_7\int_{\rho}^r\e^{\chi\tau}\tau^{3-m}\int_{B^+_{\tau}}\vert \Delta_{\gamma}u\vert \vert Du\vert \f\mu_{\gamma}\f\tau\\
&\leq 2C_7\int_{\rho}^r\e^{\chi\tau}\tau^{4-m}\int_{B^+_{\tau}}\vert \Delta_{\gamma}u\vert^2 \f\mu_{\gamma}\f\tau+2C_7\int_{\rho}^r\e^{\chi\tau}\tau^{2-m}\int_{B^+_{\tau}} \vert Du\vert^2 \f\mu_{\gamma}\f\tau\\
&\leq2C_7\int_{\rho}^r\e^{\chi\tau}\tau^{4-m}\int_{B^+_{\tau}}\vert \Delta_{\gamma}u\vert^2 \f\mu_{\gamma}\f\tau+ C_7\int_{\rho}^r\e^{\chi\tau}\tau^{4-m}\int_{B^+_{\tau}} \vert Du\vert^4 \f\mu_{\gamma}\f\tau\\
&\quad+C_7\e^{\chi R}G^{m/2}\mathcal{L}^m(B_1^+)\cdot (r-\rho) \numberthis\label{95}
\end{align*}
where we applied Young's inequality in the last step. Applying the interpolation inequality \eqref{72} yields
\begin{align*}
RHS\eqref{95}&\leq C_7C_{25}\Vert u\Vert^2_{L^{\infty}(B_1^+)}\int_{\rho}^r\e^{\chi\tau}\tau^{4-m}\int_{B^+_{\tau}} \vert D^2u\vert^2 \f\mu_{\gamma}\f\tau+C_{32}\cdot (r-\rho)\numberthis\label{96}
\end{align*}
where $C_{32}=C_{32}(\Vert Dg\Vert_{C^2},\Vert u\Vert_{L^{\infty}(B_1^+)})$. Hence,
\begin{align*}
RHS\eqref{93}&\leq 2C_7\int_{\rho}^r\e^{\chi\tau}\tau^{4-m}\int_{B^+_{\tau}}\vert \Delta_{\gamma}u\vert^2 \f\mu_{\gamma}\f\tau\\
&\quad+C_7C_{25}\Vert u\Vert^2_{L^{\infty}(B_1^+)}\int_{\rho}^r\e^{\chi\tau}\tau^{4-m}\int_{B^+_{\tau}} \vert D^2u\vert^2 \f\mu_{\gamma}\f\tau+C_{32}\cdot (r-\rho)\numberthis\label{97}.
\end{align*}
For the last integral in \eqref{93} we obtain due to $\vert \Delta_{\gamma} u\vert\leq G\vert D^2u\vert+C_7\vert Du\vert$,
\begin{align*}
&4H\int_{B^+_r\setminus B_{\rho}^+}\e^{\chi\vert x\vert}\dfrac{\vert \Delta_{\gamma}u\vert \vert D^2u\vert}{ \vert x\vert^{m-5}} \f\mu_{\gamma}\\
&\leq 8C_5\int_{B^+_r\setminus B_{\rho}^+}\e^{\chi\vert x\vert}\dfrac{\vert D^2u\vert^2 }{ \vert x\vert^{m-5}} \f\mu_{\gamma}+24C_5^2\int_{B^+_r\setminus B_{\rho}^+}\e^{\chi\vert x\vert}\dfrac{\vert Du\vert \vert D^2u\vert}{ \vert x\vert^{m-5}} \f\mu_{\gamma}\\
&\leq (8C_5+12C_5^2)\int_{B^+_r\setminus B_{\rho}^+}\e^{\chi\vert x\vert}\dfrac{\vert D^2u\vert^2 }{ \vert x\vert^{m-5}} \f\mu_{\gamma}+12C_5^2R^2\int_{B^+_r\setminus B_{\rho}^+}\e^{\chi\vert x\vert}\dfrac{\vert Du\vert^2}{ \vert x\vert^{m-3}} \f\mu_{\gamma}\numberthis\label{99}
\end{align*}
We get the second inequality in \eqref{99} due to $\vert x\vert<R$. Additonally, we rewrite the third and fourth integral in \eqref{93} with $\Delta_{\gamma}u=\Delta u+(\gamma^{ij}-\delta^{ij})u_{ij}-\gamma^{ij}\Gamma_{ij}^lu_l=:\Delta u+\Delta' u$ into
\begin{align*}
&\int_{S^+_r\setminus S_{\rho}^+}\e^{\chi\vert x\vert}\dfrac{2\Delta_{\gamma}u\cdot u_k x^k}{\vert x\vert^{m-3}} \sqrt{\gamma}\f \h^{m-1}-\int_{B^+_r\setminus B_{\rho}^+}\e^{\chi\vert x\vert}\left(\dfrac{4\Delta_{\gamma}u\cdot u_{jk}x^jx^k}{\vert x\vert^{m-2}}+\dfrac{8\Delta_{\gamma}u\cdot u_k x^k}{\vert x\vert^{m-2}}\right) \f\mu_{\gamma}\\
&=\int_{S^+_r\setminus S_{\rho}^+}\e^{\chi\vert x\vert}\dfrac{2\Delta u\cdot u_k x^k}{\vert x\vert^{m-3}} \sqrt{\gamma}\f \h^{m-1}-\int_{B^+_r\setminus B_{\rho}^+}\e^{\chi\vert x\vert}\left(\dfrac{4\Delta u\cdot u_{jk}x^jx^k}{\vert x\vert^{m-2}}+\dfrac{8\Delta u\cdot u_k x^k}{\vert x\vert^{m-2}}\right) \f\mu_{\gamma}\\
&\quad+\int_{S^+_r\setminus S_{\rho}^+}\e^{\chi\vert x\vert}\dfrac{2\Delta'u\cdot u_k x^k}{\vert x\vert^{m-3}} \sqrt{\gamma}\f \h^{m-1}-\int_{B^+_r\setminus B_{\rho}^+}\e^{\chi\vert x\vert}\left(\dfrac{4\Delta'u\cdot u_{jk}x^jx^k}{\vert x\vert^{m-2}}+\dfrac{8\Delta'u\cdot u_k x^k}{\vert x\vert^{m-2}}\right) \f\mu_{\gamma}\numberthis\label{100}.
\end{align*}

\textit{Step 4.} Next, we prove the following
\begin{lem}\label{lem2.6}
For arbitrary maps $u\in W^{2,2}(B^+)$ and $g\in C^3(B^+)$ with $(u,Du)\vert_{T_R}=(g,Dg)\vert_{T_R}$ in the sense of trace it holds
\begin{align*}
&\int_{S^+_{r}\setminus S^+_{\rho}}\e^{\chi \vert x\vert}\dfrac{2u_{ii}\cdot u_kx^k}{\vert x\vert^{m-3}}\sqrt{\gamma}\f\h^{m-1}-\int_{B^+_r\setminus B^+_{\rho}}\e^{\chi\vert x\vert}\left( \dfrac{4 u_{ii}\cdot u_{jk}x^jx^k}{\vert x\vert^{m-2}}+\dfrac{8 u_{ii}\cdot u_kx^k}{\vert x\vert^{m-2}}\right)\f\mu_{\gamma}\\
&\leq-2\int_{S^+_r\setminus S^+_{\rho}}\e^{\chi\vert x\vert}\left(-\dfrac{ u_{i} u_{ik}x^k}{\vert x\vert^{m-3}}+ 2\dfrac{ (u_{i}x^i)^2}{\vert x\vert^{m-1}}-2\dfrac{ (u_{i})^2}{\vert x\vert^{m-3}}\right)\sqrt{\gamma}\f\h^{m-1}\\
&\quad-4\int_{B^+_r\setminus B^+_{\rho}}\e^{\chi\vert x\vert}\left( \dfrac{ (u_i+u_{ij}x^j)^2}{\vert x\vert^{m-2}}+\dfrac{(m-2)(u_{i}x^i)^2}{\vert x\vert^{m}}\right)\f\mu_{\gamma}+C_{36}(m,G,\Vert Dg\Vert_{C^1},R)\cdot (r-\rho)\\
&\quad+2\chi\int_{ B^+_{r}\setminus  B^+_{\rho}}\e^{\chi \vert x\vert}\left(-\dfrac{ u_{i} u_{ik}x^k}{\vert x\vert^{m-3}}+2\dfrac{ (u_{i} x^i)^2}{\vert x\vert^{m-1}}-2\dfrac{ (u_{i})^2}{\vert x\vert^{m-3}}+\dfrac{u_{ii} u_{k}x^k}{\vert x\vert^{m-3}}\right)\f\mu_{\gamma}\\
&\quad+ 4C_5\int_{ B^+_{r}\setminus B^+_{\rho}}\e^{\chi \vert x\vert}\dfrac{\vert D^2 u\vert^2  }{\vert x\vert^{m-5}}\f\mu_{\gamma}+12C_5\int_{ B^+_{r}\setminus B^+_{\rho}}\e^{\chi \vert x\vert}\dfrac{\vert Du\vert^2  }{\vert x\vert^{m-3}}\f\mu_{\gamma}\numberthis\label{101},
\end{align*}
for a.e. $0<\rho<r<R$.
\end{lem}

To prove \Cref{lem2.6} we have to apply integration by parts. Hereby, derivatives of third order appear temporarily in intermediate steps. But $u$ is a $W^{2,2}$-map. Thus, we approximate $u$ by $C^{3}(B^+_{R-\varepsilon}) \ni u^{\varepsilon}:=g+\left[\zeta(x^m)\cdot (u-g)\right]\ast\eta_{\varepsilon}$, where $\eta_{\varepsilon}(x):=\varepsilon^{-m}\eta(\dfrac{x}{\varepsilon})$ , $\eta\in C^{\infty}(\R^m, [0,\infty))$ with $\s({\eta})\subset \overline{B}_1(0)$, $\displaystyle\int_{\R^n}\eta \f x=1$ is a mollifier and $\zeta$ is a cut-off function with $\zeta=0$ on $[0,2\varepsilon]$, $\zeta=1$ on $[3\varepsilon,\infty)$ and $\vert D\zeta\vert\leq \dfrac{2}{\varepsilon}$, $\vert D^2\zeta\vert\leq \dfrac{c}{\varepsilon^2}$. $u^{\varepsilon}$ satisfies the boundary conditions $u^{\varepsilon}=g$, $Du^{\varepsilon}=Dg$ and $D^2u^{\varepsilon}=D^2g$ on $T_{R-\varepsilon}\times [0,2\varepsilon]$. From standard properties of mollifications and Poincaré's inequality we infer $(u^{\varepsilon}-g)\rightarrow (u-g)$ in $W^{2,2}$.
\begin{proof}[Proof of \Cref{lem2.6}]
We proceed as in \cite{1} page 291 and approximate $u$ by $u^{\varepsilon}$ as already mentioned. We start with reformulation of the following boundary integral:
\begin{align*}
&\int_{S^+_{r}\setminus S^+_{\rho}}\e^{\chi \vert x\vert}\dfrac{2 u^{\varepsilon}_{ii} u^{\varepsilon}_kx^k}{\vert x\vert^{m-3}}\sqrt{\gamma}\f\h^{m-1}\\
&=\int_{\partial B^+_{r}\setminus \partial B^+_{\rho}}\e^{\chi \vert x\vert}\dfrac{2 u^{\varepsilon}_{ii} u^{\varepsilon}_kx^k}{\vert x\vert^{m-3}}\sqrt{\gamma}\f\h^{m-1}-\int_{T_{r}\setminus T_{\rho}}\e^{\chi \vert x\vert}\dfrac{2 g_{ii} g_kx^k}{\vert x\vert^{m-3}}\sqrt{\gamma}\f\h^{m-1}\\
&\leq \int_{S^+_{r}\setminus S^+_{\rho}}\e^{\chi \vert x\vert}\dfrac{2 u^{\varepsilon}_{ii} u^{\varepsilon}_kx^kx^j}{\vert x\vert^{m-2}}\nu^{j}\sqrt{\gamma}\f\h^{m-1}+2\Vert D^2g\Vert_{\infty}\Vert Dg\Vert_{\infty}G^{m/2}\int_{T_{r}\setminus T_{\rho}}\dfrac{\e^{\chi \vert x\vert}}{\vert x\vert^{m-4}}\f\h^{m-1}\\
&= \int_{B^+_{r}\setminus B^+_{\rho}}\left(\e^{\chi \vert x\vert}\dfrac{2 u^{\varepsilon}_{ii} u^{\varepsilon}_kx^kx^j}{\vert x\vert^{m-2}}\sqrt{\gamma}\right)_j\f\mathcal{L}^{m}+2\Vert D^2g\Vert_{\infty}\Vert Dg\Vert_{\infty}G^{m/2}\h^{m-1}(T_1)\left(\e^{\chi r}r^3 -\e^{\chi \rho}\rho^{3}\right)\\
&\leq  \int_{B^+_{r}\setminus B^+_{\rho}}\left(\e^{\chi \vert x\vert}\dfrac{2 u^{\varepsilon}_{ii} u^{\varepsilon}_kx^kx^j}{\vert x\vert^{m-2}}\sqrt{\gamma}\right)_j\f\mathcal{L}^{m}+C_{33}\cdot \left(r-\rho\right)\numberthis\label{102}
\end{align*}
where $C_{33}=C_{33}(\Vert Dg\Vert_{C^2},\Vert u\Vert_{L^{\infty}(B_1^+)})$. Applying product and chain rule we obtain
\begin{align*}
&\int_{ B^+_{r}\setminus  B^+_{\rho}}\left(\e^{\chi \vert x\vert}\dfrac{2u^{\varepsilon}_{ii} u^{\varepsilon}_kx^kx^j}{\vert x\vert^{m-2}}\sqrt{\gamma}\right)_j\f\mathcal{L}^m\\
&=\int_{ B^+_{r}\setminus  B^+_{\rho}}\e^{\chi \vert x\vert}\left(\dfrac{2u^{\varepsilon}_{iij}x^j u^{\varepsilon}_kx^k}{\vert x\vert^{m-2}}+\dfrac{2u^{\varepsilon}_{ii} u^{\varepsilon}_{jk}x^jx^k}{\vert x\vert^{m-2}}+\dfrac{6u^{\varepsilon}_{ii} u^{\varepsilon}_{k}x^k}{\vert x\vert^{m-2}}\right)\f\mu_{\gamma}\\
&\quad+\int_{ B^+_{r}\setminus  B^+_{\rho}}\e^{\chi \vert x\vert}\left(\dfrac{2\chi u^{\varepsilon}_{ii} u^{\varepsilon}_{k}x^k}{\vert x\vert^{m-3}}+\dfrac{2u^{\varepsilon}_{ii} u^{\varepsilon}_{k}x^k}{\vert x\vert^{m-2}}\cdot\dfrac{x^j(\sqrt{\gamma})_j}{\sqrt{\gamma}}\right)\f\mu_{\gamma}=:IX_a+IX_b \numberthis\label{103}.
\end{align*}
Hence,
\begin{align*}
&IX_a-\int_{B^+_r\setminus B^+_{\rho}}\e^{\chi\vert x\vert}\left( \dfrac{4 u^{\varepsilon}_{ii}\cdot u^{\varepsilon}_{jk}x^jx^k}{\vert x\vert^{m-2}}+\dfrac{8 u^{\varepsilon}_{ii}\cdot u^{\varepsilon}_kx^k}{\vert x\vert^{m-2}}\right)\f\mu_{\gamma}\\
&=-2\int_{ B^+_{r}\setminus  B^+_{\rho}}\e^{\chi \vert x\vert}\left(\dfrac{u^{\varepsilon}_{ii} u^{\varepsilon}_{jk}x^jx^k}{\vert x\vert^{m-2}}-\dfrac{u^{\varepsilon}_{iij}x^j u^{\varepsilon}_kx^k}{\vert x\vert^{m-2}}+\dfrac{u^{\varepsilon}_{ii} u^{\varepsilon}_{k}x^k}{\vert x\vert^{m-2}}\right)\f\mu_{\gamma}=:X\numberthis\label{104}.
\end{align*}
Now, we compute using integration by parts
\begin{align*}
&\int_{B^+_r\setminus B^+_{\rho}}\e^{\chi\vert x\vert}\left( \dfrac{ u^{\varepsilon}_{ii} u^{\varepsilon}_{jk}x^jx^k}{\vert x\vert^{m-2}}\right)\f\mu_{\gamma}\\
&=\int_{\partial B^+_r\setminus \partial B^+_{\rho}}\e^{\chi\vert x\vert}\left( \dfrac{ u^{\varepsilon}_{i} u^{\varepsilon}_{jk}x^jx^k}{\vert x\vert^{m-2}}\nu^i\right)\sqrt{\gamma}\f\h^{m-1}\\
&\quad-\int_{B^+_r\setminus B^+_{\rho}}\e^{\chi\vert x\vert}\left( \dfrac{ u^{\varepsilon}_{i} u^{\varepsilon}_{ijk}x^jx^k}{\vert x\vert^{m-2}}+\dfrac{ 2u^{\varepsilon}_{i} u^{\varepsilon}_{ik}x^k}{\vert x\vert^{m-2}}+\dfrac{ (2-m)u^{\varepsilon}_{i}x^i u^{\varepsilon}_{jk}x^jx^k}{\vert x\vert^{m}}\right)\f\mu_{\gamma}\\
&\quad-\int_{ B^+_{r}\setminus  B^+_{\rho}}\e^{\chi \vert x\vert}\left(\dfrac{\chi u^{\varepsilon}_{i}x^i u^{\varepsilon}_{jk}x^jx^k}{\vert x\vert^{m-1}}+\dfrac{u^{\varepsilon}_{jk}x^jx^k}{\vert x\vert^{m-2}}\cdot\dfrac{u^{\varepsilon}_{i}(\sqrt{\gamma})_i}{\sqrt{\gamma}}\right)\f\mu_{\gamma},\numberthis\label{105}
\end{align*}
a.e. $\rho$ and $r$, where $\nu^i$ denotes the $i$-th component of the unit normal vector. We have $\nu_S=\dfrac{x}{\vert x\vert}$ on $S^+_R$ and $\nu_T=-e_{m}=-(0,\ldots,0,1)$ on $T_R$. We split the boundary integral in \eqref{105} into flat and curved part, and keep the boundary conditions in mind. It holds
\begin{align*}
&\int_{B^+_r\setminus B^+_{\rho}}\e^{\chi\vert x\vert}\left( \dfrac{ u^{\varepsilon}_{ii} u^{\varepsilon}_{jk}x^jx^k}{\vert x\vert^{m-2}}\right)\f\mu_{\gamma}\\
&=\int_{S^+_r\setminus S^+_{\rho}}\e^{\chi\vert x\vert}\left( \dfrac{ u^{\varepsilon}_{i}x^i u^{\varepsilon}_{jk}x^jx^k}{\vert x\vert^{m-1}}\right)\sqrt{\gamma}\f\h^{m-1}-\int_{T_r\setminus T_{\rho}}\e^{\chi\vert x\vert}\left( \dfrac{ g_{m} g_{jk}x^jx^k}{\vert x\vert^{m-2}}\right)\sqrt{\gamma}\f\h^{m-1}\\
&\quad-\int_{B^+_r\setminus B^+_{\rho}}\e^{\chi\vert x\vert}\left( \dfrac{ u^{\varepsilon}_{i} u^{\varepsilon}_{ijk}x^jx^k}{\vert x\vert^{m-2}}+\dfrac{ 2u^{\varepsilon}_{i} u^{\varepsilon}_{ik}x^k}{\vert x\vert^{m-2}}+\dfrac{ (2-m)u^{\varepsilon}_{i}x^i u^{\varepsilon}_{jk}x^jx^k}{\vert x\vert^{m}}\right)\f\mu_{\gamma}\\
&\quad-\int_{ B^+_{r}\setminus  B^+_{\rho}}\e^{\chi \vert x\vert}\left(\dfrac{\chi u^{\varepsilon}_{i}x^i u^{\varepsilon}_{jk}x^jx^k}{\vert x\vert^{m-1}}+\dfrac{u^{\varepsilon}_{jk}x^jx^k}{\vert x\vert^{m-2}}\cdot\dfrac{u^{\varepsilon}_{i}(\sqrt{\gamma})_i}{\sqrt{\gamma}}\right)\f\mu_{\gamma}.\numberthis\label{106}
\end{align*}
Next, we apply integration by parts on '$\displaystyle -\int_{B^+_r\setminus B^+_{\rho}}\e^{\chi\vert x\vert}\left( \dfrac{ u^{\varepsilon}_{i} u^{\varepsilon}_{ijk}x^jx^k}{\vert x\vert^{m-2}}\right)\f\mu_{\gamma}$' in \eqref{106}:
\begin{align*}
&-\int_{B^+_r\setminus B^+_{\rho}}\e^{\chi\vert x\vert}\left( \dfrac{ u^{\varepsilon}_{i} u^{\varepsilon}_{ijk}x^jx^k}{\vert x\vert^{m-2}}\right)\f\mu_{\gamma}\\
&=-\int_{S^+_r\setminus S^+_{\rho}}\e^{\chi\vert x\vert}\left( \dfrac{ u^{\varepsilon}_{i} u^{\varepsilon}_{ik}x^k}{\vert x\vert^{m-3}}\right)\sqrt{\gamma}\f\h^{m-1}+\int_{T_r\setminus T_{\rho}}\e^{\chi\vert x\vert}\left( \dfrac{ g_{i} g_{ik}x^kx^m}{\vert x\vert^{m-2}}\right)\sqrt{\gamma}\f\h^{m-1}\\
&\quad+\int_{B^+_r\setminus B^+_{\rho}}\e^{\chi\vert x\vert}\left( \dfrac{ u^{\varepsilon}_{ij}x^j u^{\varepsilon}_{ik}x^k}{\vert x\vert^{m-2}}+3\dfrac{ u^{\varepsilon}_{i} u_{ik}x^k}{\vert x\vert^{m-2}}\right)\f\mu_{\gamma}\\
&\quad+\int_{ B^+_{r}\setminus  B^+_{\rho}}\e^{\chi \vert x\vert}\left(\dfrac{\chi u^{\varepsilon}_{i} u^{\varepsilon}_{ik}x^k}{\vert x\vert^{m-3}}+\dfrac{u^{\varepsilon}_iu^{\varepsilon}_{ik}x^k}{\vert x\vert^{m-2}}\cdot\dfrac{x^j(\sqrt{\gamma})_j}{\sqrt{\gamma}}\right)\f\mu_{\gamma}\numberthis\label{107}
\end{align*}
Using \eqref{107}, equation \eqref{106} becomes
\begin{align*}
&\int_{B^+_r\setminus B^+_{\rho}}\e^{\chi\vert x\vert}\left( \dfrac{ u^{\varepsilon}_{ii} u^{\varepsilon}_{jk}x^jx^k}{\vert x\vert^{m-2}}\right)\f\mu_{\gamma}\\
&=\int_{S^+_r\setminus S^+_{\rho}}\e^{\chi\vert x\vert}\left( \dfrac{ u^{\varepsilon}_{i}x^i u^{\varepsilon}_{jk}x^jx^k}{\vert x\vert^{m-1}}-\dfrac{ u^{\varepsilon}_{i} u^{\varepsilon}_{ik}x^k}{\vert x\vert^{m-3}}\right)\sqrt{\gamma}\f\h^{m-1}\\
&\quad+\int_{T_r\setminus T_{\rho}}\e^{\chi\vert x\vert}\left( \dfrac{ g_{i} g_{ik}x^kx^m}{\vert x\vert^{m-2}}-\dfrac{ g_{m} g_{jk}x^jx^k}{\vert x\vert^{m-2}}\right)\sqrt{\gamma}\f\h^{m-1}\\
&\quad+\int_{B^+_r\setminus B^+_{\rho}}\e^{\chi\vert x\vert}\left( \dfrac{ u^{\varepsilon}_{ij}x^j u^{\varepsilon}_{ik}x^k}{\vert x\vert^{m-2}}+\dfrac{ u^{\varepsilon}_{i} u^{\varepsilon}_{ik}x^k}{\vert x\vert^{m-2}}+\dfrac{ (m-2)u^{\varepsilon}_{i}x^i u^{\varepsilon}_{jk}x^jx^k}{\vert x\vert^{m}}\right)\f\mu_{\gamma}\\
&\quad+\int_{ B^+_{r}\setminus  B^+_{\rho}}\e^{\chi \vert x\vert}\left(\dfrac{\chi u^{\varepsilon}_{i} u^{\varepsilon}_{ik}x^k}{\vert x\vert^{m-3}}+\dfrac{u^{\varepsilon}_iu^{\varepsilon}_{ik}x^k}{\vert x\vert^{m-2}}\cdot\dfrac{x^j(\sqrt{\gamma})_j}{\sqrt{\gamma}}\right)\f\mu_{\gamma}\\
&\quad-\int_{ B^+_{r}\setminus  B^+_{\rho}}\e^{\chi \vert x\vert}\left(\dfrac{\chi u^{\varepsilon}_{i}x^i u^{\varepsilon}_{jk}x^jx^k}{\vert x\vert^{m-1}}+\dfrac{u^{\varepsilon}_{jk}x^jx^k}{\vert x\vert^{m-2}}\cdot\dfrac{u^{\varepsilon}_{i}(\sqrt{\gamma})_i}{\sqrt{\gamma}}\right)\f\mu_{\gamma}.\numberthis\label{108}
\end{align*} 
The integral over the flat part $T_r\setminus T_{\rho}$ of the boundary can be bounded from below by $ -2C_{33}\cdot (r-\rho)$. Hence, we obtain
\begin{align*}
&\int_{B^+_r\setminus B^+_{\rho}}\e^{\chi\vert x\vert}\left( \dfrac{ u^{\varepsilon}_{ii} u^{\varepsilon}_{jk}x^jx^k}{\vert x\vert^{m-2}}\right)\f\mu_{\gamma}\\
&\geq \int_{S^+_r\setminus S^+_{\rho}}\e^{\chi\vert x\vert}\left( \dfrac{ u^{\varepsilon}_{i}x^i u^{\varepsilon}_{jk}x^jx^k}{\vert x\vert^{m-1}}-\dfrac{ u^{\varepsilon}_{i} u^{\varepsilon}_{ik}x^k}{\vert x\vert^{m-3}}\right)\sqrt{\gamma}\f\h^{m-1}-2C_{33}\cdot (r-\rho)\\
&\quad+\int_{B^+_r\setminus B^+_{\rho}}\e^{\chi\vert x\vert}\left( \dfrac{ u^{\varepsilon}_{ij}x^j u^{\varepsilon}_{ik}x^k}{\vert x\vert^{m-2}}+\dfrac{ u^{\varepsilon}_{i} u^{\varepsilon}_{ik}x^k}{\vert x\vert^{m-2}}+\dfrac{ (m-2)u^{\varepsilon}_{i}x^i u^{\varepsilon}_{jk}x^jx^k}{\vert x\vert^{m}}\right)\f\mu_{\gamma}\\
&\quad+\int_{ B^+_{r}\setminus  B^+_{\rho}}\e^{\chi \vert x\vert}\left(\dfrac{\chi u^{\varepsilon}_{i} u^{\varepsilon}_{ik}x^k}{\vert x\vert^{m-3}}+\dfrac{u^{\varepsilon}_iu^{\varepsilon}_{ik}x^k}{\vert x\vert^{m-2}}\cdot\dfrac{x^j(\sqrt{\gamma})_j}{\sqrt{\gamma}}\right)\f\mu_{\gamma}\\
&\quad-\int_{ B^+_{r}\setminus  B^+_{\rho}}\e^{\chi \vert x\vert}\left(\dfrac{\chi u^{\varepsilon}_{i}x^i u^{\varepsilon}_{jk}x^jx^k}{\vert x\vert^{m-1}}+\dfrac{u^{\varepsilon}_{jk}x^jx^k}{\vert x\vert^{m-2}}\cdot\dfrac{u^{\varepsilon}_{i}(\sqrt{\gamma})_i}{\sqrt{\gamma}}\right)\f\mu_{\gamma}.\numberthis\label{110}
\end{align*} 
We continue as follows using again integration by parts,
\begin{align*}
&\int_{B^+_r\setminus B^+_{\rho}}\e^{\chi\vert x\vert}\left( -\dfrac{u^{\varepsilon}_{iij}x^j u^{\varepsilon}_kx^k}{\vert x\vert^{m-2}}\right)\f\mu_{\gamma}\\
&=\int_{S^+_r\setminus S^+_{\rho}}\e^{\chi\vert x\vert}\left( -\dfrac{u^{\varepsilon}_{ij}x^ix^j u^{\varepsilon}_kx^k}{\vert x\vert^{m-1}}\right)\sqrt{\gamma}\f\h^{m-1}+\int_{T_r\setminus T_{\rho}}\e^{\chi\vert x\vert}\left( \dfrac{g_{mj}x^j g_kx^k}{\vert x\vert^{m-2}}\right)\sqrt{\gamma}\f\h^{m-1}\\
&\quad+\int_{B^+_r\setminus B^+_{\rho}}\e^{\chi\vert x\vert}\left( \dfrac{u^{\varepsilon}_{ij}x^j u^{\varepsilon}_{ik}x^k}{\vert x\vert^{m-2}}+\dfrac{u^{\varepsilon}_{ii} u^{\varepsilon}_{k}x^k}{\vert x\vert^{m-2}}+\dfrac{u^{\varepsilon}_{i} u^{\varepsilon}_{ij}x^j}{\vert x\vert^{m-2}}-\dfrac{(m-2)u^{\varepsilon}_{ij}x^ix^j u^{\varepsilon}_kx^k}{\vert x\vert^{m}}\right)\f\mu_{\gamma}\\
&\quad+\int_{ B^+_{r}\setminus  B^+_{\rho}}\e^{\chi \vert x\vert}\left(\dfrac{\chi u^{\varepsilon}_{ij}x^ix^j u^{\varepsilon}_{k}x^k}{\vert x\vert^{m-1}}+\dfrac{u^{\varepsilon}_{ij}x^ju^{\varepsilon}_{k}x^k}{\vert x\vert^{m-2}}\cdot\dfrac{(\sqrt{\gamma})_i}{\sqrt{\gamma}}\right)\f\mu_{\gamma}\\
&\geq \int_{S^+_r\setminus S^+_{\rho}}\e^{\chi\vert x\vert}\left( -\dfrac{u^{\varepsilon}_{ij}x^ix^j u^{\varepsilon}_kx^k}{\vert x\vert^{m-1}}\right)\sqrt{\gamma}\f\h^{m-1}-C_{33}\cdot (r-\rho)\\
&\quad+\int_{B^+_r\setminus B^+_{\rho}}\e^{\chi\vert x\vert}\left( \dfrac{u^{\varepsilon}_{ij}x^j u^{\varepsilon}_{ik}x^k}{\vert x\vert^{m-2}}+\dfrac{u^{\varepsilon}_{ii} u^{\varepsilon}_{k}x^k}{\vert x\vert^{m-2}}+\dfrac{u^{\varepsilon}_{i} u^{\varepsilon}_{ij}x^j}{\vert x\vert^{m-2}}-\dfrac{(m-2)u^{\varepsilon}_{ij}x^ix^j u^{\varepsilon}_kx^k}{\vert x\vert^{m}}\right)\f\mu_{\gamma}\\
&\quad+\int_{ B^+_{r}\setminus  B^+_{\rho}}\e^{\chi \vert x\vert}\left(\dfrac{\chi u^{\varepsilon}_{ij}x^ix^j u^{\varepsilon}_{k}x^k}{\vert x\vert^{m-1}}+\dfrac{u^{\varepsilon}_{ij}x^ju^{\varepsilon}_{k}x^k}{\vert x\vert^{m-2}}\cdot\dfrac{(\sqrt{\gamma})_i}{\sqrt{\gamma}}\right)\f\mu_{\gamma}.\numberthis\label{113}
\end{align*}
Thus, we have
\begin{align*}
-\dfrac{1}{2}X&=\int_{B^+_r\setminus B^+_{\rho}}\e^{\chi\vert x\vert}\left( \dfrac{ u^{\varepsilon}_{ii} u^{\varepsilon}_{jk}x^jx^k}{\vert x\vert^{m-2}}-\dfrac{u^{\varepsilon}_{iij}x^j u^{\varepsilon}_kx^k}{\vert x\vert^{m-2}}+\dfrac{u^{\varepsilon}_{ii}u^{\varepsilon}_kx^k}{\vert x\vert^{m-2}}\right)\f\mu_{\gamma}\\
&\geq \int_{S^+_r\setminus S^+_{\rho}}\e^{\chi\vert x\vert}\left( -\dfrac{ u^{\varepsilon}_{i} u^{\varepsilon}_{ik}x^k}{\vert x\vert^{m-3}}\right)\sqrt{\gamma}\f\h^{m-1}-3C_{33}\cdot (r-\rho)\\
&\quad+2\int_{B^+_r\setminus B^+_{\rho}}\e^{\chi\vert x\vert}\left( \dfrac{ u^{\varepsilon}_{ij}x^j u^{\varepsilon}_{ik}x^k}{\vert x\vert^{m-2}}+\dfrac{ u^{\varepsilon}_{i} u^{\varepsilon}_{ik}x^k}{\vert x\vert^{m-2}}+\dfrac{u^{\varepsilon}_{ii}u^{\varepsilon}_kx^k}{\vert x\vert^{m-2}}\right)\f\mu_{\gamma}\\
&\quad+\int_{ B^+_{r}\setminus  B^+_{\rho}}\e^{\chi \vert x\vert}\left(\dfrac{\chi u^{\varepsilon}_{i} u^{\varepsilon}_{ik}x^k}{\vert x\vert^{m-3}}+\dfrac{u^{\varepsilon}_iu^{\varepsilon}_{ik}x^k}{\vert x\vert^{m-2}}\cdot\dfrac{x^j(\sqrt{\gamma})_j}{\sqrt{\gamma}}\right)\f\mu_{\gamma}\\
&\quad+\int_{ B^+_{r}\setminus  B^+_{\rho}}\e^{\chi \vert x\vert}\left(\dfrac{u^{\varepsilon}_{ij}x^ju^{\varepsilon}_{k}x^k}{\vert x\vert^{m-2}}\cdot\dfrac{(\sqrt{\gamma})_i}{\sqrt{\gamma}}-\dfrac{u^{\varepsilon}_{jk}x^jx^k}{\vert x\vert^{m-2}}\cdot\dfrac{u^{\varepsilon}_{i}(\sqrt{\gamma})_i}{\sqrt{\gamma}}\right)\f\mu_{\gamma}.\numberthis\label{114}
\end{align*}
Moreover, it holds
\begin{align*}
&\int_{ B^+_r\setminus  B^+_{\rho}}\e^{\chi\vert x\vert}\left( \dfrac{u^{\varepsilon}_{ii}u^{\varepsilon}_kx^k}{\vert x\vert^{m-2}}\right)\f\mu_{\gamma}\\
&=\int_{S^+_r\setminus S^+_{\rho}}\e^{\chi\vert x\vert}\left( \dfrac{ (u^{\varepsilon}_{i}x^i)^2}{\vert x\vert^{m-1}}\right)\sqrt{\gamma}\f\h^{m-1}-\int_{T_r\setminus T_{\rho}}\e^{\chi\vert x\vert}\left( \dfrac{ g_mg_{k}x^k}{\vert x\vert^{m-2}}\right)\sqrt{\gamma}\f\h^{m-1}\\
&\quad-\int_{B^+_r\setminus B^+_{\rho}}\e^{\chi\vert x\vert}\left( \dfrac{u^{\varepsilon}_{i}u^{\varepsilon}_{ik}x^k}{\vert x\vert^{m-2}}+\dfrac{(u^{\varepsilon}_{i})^2}{\vert x\vert^{m-2}}-\dfrac{(m-2)(u^{\varepsilon}_{i}x^i)^2}{\vert x\vert^{m}}\right)\f\mu_{\gamma}\\
&\quad-\int_{ B^+_{r}\setminus  B^+_{\rho}}\e^{\chi \vert x\vert}\left(\dfrac{\chi (u^{\varepsilon}_{i} x^i)^2}{\vert x\vert^{m-1}}+\dfrac{u^{\varepsilon}_iu^{\varepsilon}_{k}x^k}{\vert x\vert^{m-2}}\cdot\dfrac{(\sqrt{\gamma})_i}{\sqrt{\gamma}}\right)\f\mu_{\gamma}\numberthis\label{115}.
\end{align*}
The integral over the flat part $T_r\setminus T_{\rho}$ of the boundary can be bounded from below by $ -2C_{34}\cdot (r-\rho)$ where $C_{34}=C_{34}(\Vert Dg\Vert_{C^2},\Vert u\Vert_{L^{\infty}(B_1^+)})$. Thereby and with \eqref{115} follow
\begin{align*}
-\dfrac{1}{2}X\geq RHS\eqref{114}&\geq\int_{S^+_r\setminus S^+_{\rho}}\e^{\chi\vert x\vert}\left( 2\dfrac{ (u^{\varepsilon}_{i}x^i)^2}{\vert x\vert^{m-1}}-\dfrac{ u^{\varepsilon}_{i} u^{\varepsilon}_{ik}x^k}{\vert x\vert^{m-3}}\right)\sqrt{\gamma}\f\h^{m-1}-C_{35}\cdot (r-\rho)\\
&\quad+2\int_{B^+_r\setminus B^+_{\rho}}\e^{\chi\vert x\vert}\left( \dfrac{ (u^{\varepsilon}_{ij}x^j)^2}{\vert x\vert^{m-2}}-\dfrac{(u^{\varepsilon}_{i})^2}{\vert x\vert^{m-2}}+\dfrac{(m-2)(u^{\varepsilon}_{i}x^i)^2}{\vert x\vert^{m}}\right)\f\mu_{\gamma}\\
&\quad+\int_{ B^+_{r}\setminus  B^+_{\rho}}\e^{\chi \vert x\vert}\left(\dfrac{\chi u^{\varepsilon}_{i} u^{\varepsilon}_{ik}x^k}{\vert x\vert^{m-3}}+\dfrac{u^{\varepsilon}_iu^{\varepsilon}_{ik}x^k}{\vert x\vert^{m-2}}\cdot\dfrac{x^j(\sqrt{\gamma})_j}{\sqrt{\gamma}}\right)\f\mu_{\gamma}\\
&\quad+\int_{ B^+_{r}\setminus  B^+_{\rho}}\e^{\chi \vert x\vert}\left(\dfrac{u^{\varepsilon}_{ij}x^ju^{\varepsilon}_{k}x^k}{\vert x\vert^{m-2}}\cdot\dfrac{(\sqrt{\gamma})_i}{\sqrt{\gamma}}-\dfrac{u^{\varepsilon}_{jk}x^jx^k}{\vert x\vert^{m-2}}\cdot\dfrac{u^{\varepsilon}_{i}(\sqrt{\gamma})_i}{\sqrt{\gamma}}\right)\f\mu_{\gamma}\\
&\quad-2\int_{ B^+_{r}\setminus  B^+_{\rho}}\e^{\chi \vert x\vert}\left(\dfrac{\chi (u^{\varepsilon}_{i} x^i)^2}{\vert x\vert^{m-1}}+\dfrac{u^{\varepsilon}_iu^{\varepsilon}_{k}x^k}{\vert x\vert^{m-2}}\cdot\dfrac{(\sqrt{\gamma})_i}{\sqrt{\gamma}}\right)\f\mu_{\gamma},\numberthis\label{117}
\end{align*}
where $C_{35}=C_{35}(\Vert Dg\Vert_{C^2},\Vert u\Vert_{L^{\infty}(B_1^+)}):=3C_{33}+C_{34}$. In addition, we get by Gauss's integration theorem (cf. \cite{1}, page 292)
\begin{align*}
0&=-\int_{S^+_r\setminus S^+_{\rho}}\e^{\chi\vert x\vert}\left( 2\dfrac{ (u^{\varepsilon}_{i})^2}{\vert x\vert^{m-3}}\right)\sqrt{\gamma}\f\h^{m-1}\\
&\quad+2\int_{B^+_r\setminus B^+_{\rho}}\e^{\chi\vert x\vert}\left( 2\dfrac{ u^{\varepsilon}_iu^{\varepsilon}_{ij}x^j}{\vert x\vert^{m-2}}+2\dfrac{(u^{\varepsilon}_{i})^2}{\vert x\vert^{m-2}}+\dfrac{\chi (u^{\varepsilon}_{i})^2}{\vert x\vert^{m-3}}+\dfrac{(u^{\varepsilon}_i)^2}{\vert x\vert^{m-2}}\cdot\dfrac{x^j(\sqrt{\gamma})_j}{\sqrt{\gamma}}\right)\f\mu_{\gamma}\numberthis\label{118},
\end{align*}
which we add to \eqref{117}. Hence, after suitbale reformulations we obtain the following inequality
\begin{align*}
X&\leq-2\int_{S^+_r\setminus S^+_{\rho}}\e^{\chi\vert x\vert}\left(-\dfrac{ u^{\varepsilon}_{i} u^{\varepsilon}_{ik}x^k}{\vert x\vert^{m-3}}+ 2\dfrac{ (u^{\varepsilon}_{i}x^i)^2}{\vert x\vert^{m-1}}-2\dfrac{ (u^{\varepsilon}_{i})^2}{\vert x\vert^{m-3}}\right)\sqrt{\gamma}\f\h^{m-1}+2C_{35}\cdot (r-\rho)\\
&\quad-4\int_{B^+_r\setminus B^+_{\rho}}\e^{\chi\vert x\vert}\left( \dfrac{ (u^{\varepsilon}_i+u^{\varepsilon}_{ij}x^j)^2}{\vert x\vert^{m-2}}+\dfrac{(m-2)(u^{\varepsilon}_{i}x^i)^2}{\vert x\vert^{m}}\right)\f\mu_{\gamma}\\
&\quad+2\chi\int_{ B^+_{r}\setminus  B^+_{\rho}}\e^{\chi \vert x\vert}\left(-\dfrac{ u^{\varepsilon}_{i} u^{\varepsilon}_{ik}x^k}{\vert x\vert^{m-3}}+2\dfrac{ (u^{\varepsilon}_{i} x^i)^2}{\vert x\vert^{m-1}}-2\dfrac{ (u^{\varepsilon}_{i})^2}{\vert x\vert^{m-3}}\right)\f\mu_{\gamma}\\
&\quad-2\int_{ B^+_{r}\setminus  B^+_{\rho}}\e^{\chi \vert x\vert}\left(\dfrac{u^{\varepsilon}_{ij}x^ju^{\varepsilon}_{k}x^k}{\vert x\vert^{m-2}}\cdot\dfrac{(\sqrt{\gamma})_i}{\sqrt{\gamma}}-\dfrac{u^{\varepsilon}_{jk}x^jx^k}{\vert x\vert^{m-2}}\cdot\dfrac{u^{\varepsilon}_{i}(\sqrt{\gamma})_i}{\sqrt{\gamma}}+2\dfrac{(u^{\varepsilon}_i)^2}{\vert x\vert^{m-2}}\cdot\dfrac{x^j(\sqrt{\gamma})_j}{\sqrt{\gamma}}\right)\f\mu_{\gamma}\\
&\quad+2\int_{ B^+_{r}\setminus  B^+_{\rho}}\e^{\chi \vert x\vert}\left(-\dfrac{u^{\varepsilon}_iu^{\varepsilon}_{ik}x^k}{\vert x\vert^{m-2}}\cdot\dfrac{x^j(\sqrt{\gamma})_j}{\sqrt{\gamma}}+2\dfrac{u^{\varepsilon}_iu^{\varepsilon}_{k}x^k}{\vert x\vert^{m-2}}\cdot\dfrac{(\sqrt{\gamma})_i}{\sqrt{\gamma}}\right)\f\mu_{\gamma}.\numberthis\label{119}
\end{align*}
Thus, from standard properties of mollification we get 
\begin{align*}
&\lim_{\varepsilon\searrow0}\left(IX_b+X\right)+C_{33}\cdot (r-\rho)\\
&=\int_{S^+_{r}\setminus S^+_{\rho}}\e^{\chi \vert x\vert}\dfrac{2u_{ii}\cdot u_kx^k}{\vert x\vert^{m-3}}\sqrt{\gamma}\f\h^{m-1}-\int_{B^+_r\setminus B^+_{\rho}}\e^{\chi\vert x\vert}\left( \dfrac{4 u_{ii}\cdot u_{jk}x^jx^k}{\vert x\vert^{m-2}}+\dfrac{8 u_{ii}\cdot u_kx^k}{\vert x\vert^{m-2}}\right)\f\mu_{\gamma}\\
&\leq-2\int_{S^+_r\setminus S^+_{\rho}}\e^{\chi\vert x\vert}\left(-\dfrac{ u_{i} u_{ik}x^k}{\vert x\vert^{m-3}}+ 2\dfrac{ (u_{i}x^i)^2}{\vert x\vert^{m-1}}-2\dfrac{ (u_{i})^2}{\vert x\vert^{m-3}}\right)\sqrt{\gamma}\f\h^{m-1}+C_{36}\cdot (r-\rho)\\
&\quad-4\int_{B^+_r\setminus B^+_{\rho}}\e^{\chi\vert x\vert}\left( \dfrac{ (u_i+u_{ij}x^j)^2}{\vert x\vert^{m-2}}+\dfrac{(m-2)(u_{i}x^i)^2}{\vert x\vert^{m}}\right)\f\mu_{\gamma}\\
&\quad+2\chi\int_{ B^+_{r}\setminus  B^+_{\rho}}\e^{\chi \vert x\vert}\left(-\dfrac{ u_{i} u_{ik}x^k}{\vert x\vert^{m-3}}+2\dfrac{ (u_{i} x^i)^2}{\vert x\vert^{m-1}}-2\dfrac{ (u_{i})^2}{\vert x\vert^{m-3}}+\dfrac{u_{ii} u_{k}x^k}{\vert x\vert^{m-3}}\right)\f\mu_{\gamma}\\
&\quad-2\int_{ B^+_{r}\setminus  B^+_{\rho}}\e^{\chi \vert x\vert}\left(\dfrac{u_{ij}x^ju_{k}x^k}{\vert x\vert^{m-2}}\cdot\dfrac{(\sqrt{\gamma})_i}{\sqrt{\gamma}}-\dfrac{u_{jk}x^jx^k}{\vert x\vert^{m-2}}\cdot\dfrac{u_{i}(\sqrt{\gamma})_i}{\sqrt{\gamma}}\right)\f\mu_{\gamma}\\
&\quad-2\int_{ B^+_{r}\setminus  B^+_{\rho}}\e^{\chi \vert x\vert}\left(\dfrac{u_iu_{ik}x^k}{\vert x\vert^{m-2}}\cdot\dfrac{x^j(\sqrt{\gamma})_j}{\sqrt{\gamma}}+2\dfrac{(u_i)^2}{\vert x\vert^{m-2}}\cdot\dfrac{x^j(\sqrt{\gamma})_j}{\sqrt{\gamma}}\right)\f\mu_{\gamma}\\
&\quad+2\int_{ B^+_{r}\setminus  B^+_{\rho}}\e^{\chi \vert x\vert}\left(2\dfrac{u_iu_{k}x^k}{\vert x\vert^{m-2}}\cdot\dfrac{(\sqrt{\gamma})_i}{\sqrt{\gamma}}+\dfrac{u_{ii} u_{k}x^k}{\vert x\vert^{m-2}}\cdot\dfrac{x^j(\sqrt{\gamma})_j}{\sqrt{\gamma}}\right)\f\mu_{\gamma} \numberthis\label{120}
\end{align*}
for a.e. $\rho$ and $r$ where $C_{36}=C_{36}(\Vert Dg\Vert_{C^2},\Vert u\Vert_{L^{\infty}(B_1^+)}):=2C_{35}+C_{33}$. We estimate the last three integrals in \eqref{120} as follows,
\begin{align*}
&\quad-2\int_{ B^+_{r}\setminus  B^+_{\rho}}\e^{\chi \vert x\vert}\left(\dfrac{u_{ij}x^ju_{k}x^k}{\vert x\vert^{m-2}}\cdot\dfrac{(\sqrt{\gamma})_i}{\sqrt{\gamma}}-\dfrac{u_{jk}x^jx^k}{\vert x\vert^{m-2}}\cdot\dfrac{u_{i}(\sqrt{\gamma})_i}{\sqrt{\gamma}}\right)\f\mu_{\gamma}\\
&\quad-2\int_{ B^+_{r}\setminus  B^+_{\rho}}\e^{\chi \vert x\vert}\left(\dfrac{u_iu_{ik}x^k}{\vert x\vert^{m-2}}\cdot\dfrac{x^j(\sqrt{\gamma})_j}{\sqrt{\gamma}}+2\dfrac{(u_i)^2}{\vert x\vert^{m-2}}\cdot\dfrac{x^j(\sqrt{\gamma})_j}{\sqrt{\gamma}}\right)\f\mu_{\gamma}\\
&\quad+2\int_{ B^+_{r}\setminus  B^+_{\rho}}\e^{\chi \vert x\vert}\left(2\dfrac{u_iu_{k}x^k}{\vert x\vert^{m-2}}\cdot\dfrac{(\sqrt{\gamma})_i}{\sqrt{\gamma}}+\dfrac{u_{ii} u_{k}x^k}{\vert x\vert^{m-2}}\cdot\dfrac{x^j(\sqrt{\gamma})_j}{\sqrt{\gamma}}\right)\f\mu_{\gamma} \\
&\leq 8C_5\int_{ B^+_{r}\setminus B_{\rho}^+}\e^{\chi \vert x\vert}\left(\dfrac{\vert D^2 u\vert \vert Du\vert  }{\vert x\vert^{m-4}}+\dfrac{ \vert Du\vert^2  }{\vert x\vert^{m-3}}\right)\f\mu_{\gamma}\\
&\leq 4C_5\int_{ B^+_{r}\setminus B_{\rho}^+}\e^{\chi \vert x\vert}\dfrac{\vert D^2 u\vert^2  }{\vert x\vert^{m-5}}\f\mu_{\gamma}+12C_5\int_{ B^+_{r}\setminus B_{\rho}^+}\e^{\chi \vert x\vert}\dfrac{\vert Du\vert^2  }{\vert x\vert^{m-3}}\f\mu_{\gamma}\numberthis\label{121}.
\end{align*}
Altogether, we have
\begin{align*}
&\int_{S^+_{r}\setminus S^+_{\rho}}\e^{\chi \vert x\vert}\dfrac{2u_{ii}\cdot u_kx^k}{\vert x\vert^{m-3}}\sqrt{\gamma}\f\h^{m-1}-\int_{B^+_r\setminus B^+_{\rho}}\e^{\chi\vert x\vert}\left( \dfrac{4 u_{ii}\cdot u_{jk}x^jx^k}{\vert x\vert^{m-2}}+\dfrac{8 u_{ii}\cdot u_kx^k}{\vert x\vert^{m-2}}\right)\f\mu_{\gamma}\\
&\leq-2\int_{S^+_r\setminus S^+_{\rho}}\e^{\chi\vert x\vert}\left(-\dfrac{ u_{i} u_{ik}x^k}{\vert x\vert^{m-3}}+ 2\dfrac{ (u_{i}x^i)^2}{\vert x\vert^{m-1}}-2\dfrac{ (u_{i})^2}{\vert x\vert^{m-3}}\right)\sqrt{\gamma}\f\h^{m-1}+C_{36}\cdot (r-\rho)\\
&\quad-4\int_{B^+_r\setminus B^+_{\rho}}\e^{\chi\vert x\vert}\left( \dfrac{ (u_i+u_{ij}x^j)^2}{\vert x\vert^{m-2}}+\dfrac{(m-2)(u_{i}x^i)^2}{\vert x\vert^{m}}\right)\f\mu_{\gamma}\\
&\quad+2\chi\int_{ B^+_{r}\setminus  B^+_{\rho}}\e^{\chi \vert x\vert}\left(-\dfrac{ u_{i} u_{ik}x^k}{\vert x\vert^{m-3}}+2\dfrac{ (u_{i} x^i)^2}{\vert x\vert^{m-1}}-2\dfrac{ (u_{i})^2}{\vert x\vert^{m-3}}+\dfrac{u_{ii} u_{k}x^k}{\vert x\vert^{m-3}}\right)\f\mu_{\gamma}\\
&\quad +4C_5\int_{ B^+_{r}\setminus B_{\rho}^+}\e^{\chi \vert x\vert}\dfrac{\vert D^2 u\vert^2  }{\vert x\vert^{m-5}}\f\mu_{\gamma}+12C_5\int_{ B^+_{r}\setminus B_{\rho}^+}\e^{\chi \vert x\vert}\dfrac{\vert Du\vert^2  }{\vert x\vert^{m-3}}\f\mu_{\gamma}\numberthis\label{122}.
\end{align*}
This concludes the proof of \Cref{lem2.6}.
\end{proof}
\textit{Step 5}. We continue with an estimate of the last integral in \eqref{100}. Since $\vert \Delta'u\vert \leq H\vert x\vert \vert D^2u\vert+C_7\vert Du \vert$, it holds
\begin{align*}
&-\int_{B^+_r\setminus B_{\rho}^+}\e^{\chi\vert x\vert}\left(\dfrac{4\Delta'u\cdot u_{jk}x^jx^k}{\vert x\vert^{m-2}}+\dfrac{8\Delta'u\cdot u_k x^k}{\vert x\vert^{m-2}}\right) \f\mu_{\gamma}\\
&\leq \int_{B^+_r\setminus B_{\rho}^+}\e^{\chi \vert x\vert}\left(\dfrac{4(H\vert x\vert \vert D^2u\vert +C_7\vert Du\vert ) \vert D^2u\vert }{\vert x\vert^{m-4}}+\dfrac{8(H\vert x\vert \vert D^2u\vert +C_7\vert Du\vert ) \vert Du\vert }{\vert x\vert^{m-3}}\right)\f\mu_{\gamma}\\
&\leq C_{37}\int_{B^+_r\setminus B_{\rho}^+}\e^{\chi \vert x\vert}\dfrac{ \vert D^2u\vert^2 }{\vert x\vert^{m-5}}\f\mu_{\gamma}+C_{38}\int_{B_{r}^+\setminus B_{\rho}^+}\e^{\chi \vert x\vert}\dfrac{\vert Du\vert^2}{\vert x\vert^{m-3}} \f\mu_{\gamma}\numberthis\label{124}
\end{align*}
where $C_{37}:=8H+2C_7$ and $C_{38}:=4H+6C_7$. Furthermore, we estimate the second to last integral in \eqref{100}. Due to
\begin{align*}
\vert \Delta' u\vert=\vert(\gamma^{ij}-\delta^{ij})u_{ij}-\gamma^{ij}\Gamma_{ij}^lu_l\vert\leq H\vert x\vert \vert D^2u\vert+C_7\vert Du\vert\numberthis\label{125}
\end{align*}
we get the following estimate,
\begin{align*}
&\int_{S^+_r\setminus S_{\rho}^+}\e^{\chi\vert x\vert}\dfrac{2\Delta'u\cdot u_k x^k}{\vert x\vert^{m-3}} \sqrt{\gamma}\f \h^{m-1}\leq \int_{S_{r}^+\cup S_{\rho}^+}\e^{\chi \vert x\vert}\dfrac{2(H\vert x\vert \vert D^2u\vert +C_7\vert Du\vert ) \vert Du\vert }{\vert x\vert^{m-4}} \sqrt{\gamma}\f \h^{m-1}\\
&\leq H\int_{S_{r}^+\cup S_{\rho}^+}\e^{\chi \vert x\vert}\dfrac{\vert D^2u\vert^2  }{\vert x\vert^{m-6}} \sqrt{\gamma}\f \h^{m-1}+(2C_7+H)\int_{S_{r}^+\cup S_{\rho}^+}\e^{\chi \vert x\vert}\dfrac{\vert Du\vert^2}{\vert x\vert^{m-4}} \sqrt{\gamma}\f \h^{m-1}\numberthis\label{126},
\end{align*}
where we used $2\vert D^2u\vert \vert Du\vert\leq \vert D^2u\vert^2\vert x\vert+ \vert Du\vert^2/\vert x\vert$ in the last step. Further, we rewrite the second integral in \eqref{93} with $\Delta_{\gamma}=\Delta+\Delta'$ as follows,
\begin{align*}
-2\chi\int_{B^+_r\setminus B_{\rho}^+}\e^{\chi\vert x\vert}\dfrac{\Delta_{\gamma}u\cdot u_k x^k}{\vert x\vert^{m-3}}&=-2\chi\int_{B^+_r\setminus B_{\rho}^+}\e^{\chi\vert x\vert}\dfrac{u_{ii} u_k x^k}{\vert x\vert^{m-3}} \f\mu_{\gamma}\underbrace{-2\chi\int_{B^+_r\setminus B_{\rho}^+}\e^{\chi\vert x\vert}\dfrac{\Delta'u\cdot u_k x^k}{\vert x\vert^{m-3}}\f\mu_{\gamma}}_{=:XI}\numberthis\label{127}.
\end{align*}
It follows due to \eqref{125}:
\begin{align*}
 XI &\leq 2H\chi\int_{B^+_r\setminus B_{\rho}^+}\e^{\chi \vert x\vert }\dfrac{\vert D^2u\vert \vert Du\vert }{\vert x\vert^{m-5}} \f\mu_{\gamma}+2C_7\chi\int_{B^+_r\setminus B_{\rho}^+}\e^{\chi \vert x\vert }\dfrac{\vert Du\vert^2 }{\vert x\vert^{m-4}} \f\mu_{\gamma}\numberthis\label{128}
\end{align*}
Using $2\vert D^2u\vert \vert Du\vert\leq \vert D^2u\vert^2+ \vert Du\vert^2$ yields
\begin{align*}
RHS\eqref{128}&\leq H\chi\int_{B^+_r\setminus B_{\rho}^+}\e^{\chi \vert x\vert }\dfrac{\vert D^2u\vert^2 }{\vert x\vert^{m-5}} \f\mu_{\gamma}+H\chi\int_{B^+_r\setminus B_{\rho}^+}\e^{\chi \vert x\vert }\dfrac{ \vert Du\vert^2 }{\vert x\vert^{m-5}} \f\mu_{\gamma}\\
&\quad+2C_7\chi\int_{B^+_r\setminus B_{\rho}^+}\e^{\chi \vert x\vert }\dfrac{\vert Du\vert^2 }{\vert x\vert^{m-4}} \f\mu_{\gamma}\\
&\leq H\chi\int_{B^+_r\setminus B_{\rho}^+}\e^{\chi \vert x\vert }\dfrac{\vert D^2u\vert^2 }{\vert x\vert^{m-5}} \f\mu_{\gamma}+(2C_7R+HR^2)\chi\int_{B^+_r\setminus B_{\rho}^+}\e^{\chi \vert x\vert }\dfrac{ \vert Du\vert^2 }{\vert x\vert^{m-3}} \f\mu_{\gamma}\numberthis\label{130}.
\end{align*}
Because of \eqref{97}, \eqref{99}, \eqref{101}, \eqref{124}, \eqref{126}, \eqref{127} and \eqref{130} the following estimate holds,
\begin{align*}
&\lim_{\nu\rightarrow\infty}\int_{\rho}^r\e^{\chi\tau}\tau^{3-m}VII\f\tau\leq  RHS\eqref{93}\\
&\leq \int_{\rho}^r\e^{\chi\tau}I'(\tau)\f\tau\bcancel{-2\chi\int_{B^+_r\setminus B_{\rho}^+}\e^{\chi\vert x\vert}\dfrac{u_{ii} u_k x^k}{\vert x\vert^{m-3}} \f\mu_{\gamma}}+C_{39}\cdot (r-\rho)\\
&\quad-2\int_{S^+_r\setminus S^+_{\rho}}\e^{\chi\vert x\vert}\left(-\dfrac{ u_{i} u_{ik}x^k}{\vert x\vert^{m-3}}+ 2\dfrac{ (u_{i}x^i)^2}{\vert x\vert^{m-1}}-2\dfrac{ (u_{i})^2}{\vert x\vert^{m-3}}\right)\sqrt{\gamma}\f\h^{m-1}\\
&\quad-4\int_{B^+_r\setminus B^+_{\rho}}\e^{\chi\vert x\vert}\left( \dfrac{ (u_i+u_{ij}x^j)^2}{\vert x\vert^{m-2}}+\dfrac{(m-2)(u_{i}x^i)^2}{\vert x\vert^{m}}\right)\f\mu_{\gamma}\\
&\quad+2\chi\int_{ B^+_{r}\setminus  B^+_{\rho}}\e^{\chi \vert x\vert}\left(-\dfrac{ u_{i} u_{ik}x^k}{\vert x\vert^{m-3}}+2\dfrac{ (u_{i} x^i)^2}{\vert x\vert^{m-1}}-2\dfrac{ (u_{i})^2}{\vert x\vert^{m-3}}+\bcancel{\dfrac{u_{ii} u_{k}x^k}{\vert x\vert^{m-3}}}\right)\f\mu_{\gamma}\\
&\quad +C_{40}\int_{ B^+_{r}\setminus B_{\rho}^+}\e^{\chi\vert x\vert}\dfrac{\vert D^2 u\vert^2  }{\vert x\vert^{m-5}}\f\mu_{\gamma}+C_{41}\int_{ B^+_{r}\setminus B_{\rho}^+}\e^{\chi\vert x\vert}\dfrac{\vert Du\vert^2  }{\vert x\vert^{m-3}}\f\mu_{\gamma}\\
&\quad+H\int_{S_{r}^+\cup S_{\rho}^+}\e^{\chi\vert x\vert}\dfrac{\vert D^2u\vert^2  }{\vert x\vert^{m-6}} \sqrt{\gamma}\f\h^{m-1}+(2C_7+H)\int_{S_{r}^+\cup S_{\rho}^+}\e^{\chi\vert x\vert}\dfrac{\vert Du\vert^2}{\vert x\vert^{m-4}}\sqrt{\gamma} \f\h^{m-1}\\
&\quad+2C_7\int_{\rho}^r\e^{\chi\tau}\tau^{4-m}\int_{B^+_{\tau}}\vert \Delta_{\gamma}u\vert^2 \f\mu_{\gamma}\f\tau+C_{42}\int_{\rho}^r\e^{\chi\tau}\tau^{4-m}\int_{B^+_{\tau}} \vert D^2u\vert^2 \f\mu_{\gamma}\f\tau\numberthis\label{131}
\end{align*}
where $C_{39}=C_{39}(\Vert Dg\Vert_{C^2},\Vert u\Vert_{L^{\infty}(B_1^+)}):=C_{32}+C_{36}$, $C_{40}=C_{40}(\Vert Dg\Vert_{C^2},\Vert u\Vert_{L^{\infty}(B_1^+)}):=8C_5+12C_5^2+H\chi+C_{37}$, $C_{41}=C_{41}(\Vert Dg\Vert_{C^2},\Vert u\Vert_{L^{\infty}(B_1^+)}):=(2C_7R+HR^2)\chi+12C_5^2R^2+C_{38}$ and $C_{42}=C_{42}(\Vert u\Vert_{L^{\infty}(B_1^+)}):=C_7C_{25}\Vert u\Vert^2_{L^{\infty}(B_1^+)}$.

\textit{Step 6}. Observe that it holds
\begin{align*}
\xi_{ij}^k=\dfrac{x^ix^jx^k}{\tau^2\vert x\vert^2}\si''+\dfrac{1}{\tau}\left(\delta_{jk}x^i+\delta_{ik}x^j+\delta_{ij}x^k-\dfrac{x^ix^jx^k}{\vert x\vert^2}\right)\si'\cdot\dfrac{1}{\vert x\vert}\numberthis\label{132}
\end{align*}
where $\si'':=\psi_{\nu}''\left(\dfrac{\vert x\vert}{\tau}\right)$. Putting \eqref{132} into $VIII$ yields
\begin{align*}
VIII&=2\int_{B^+}\Delta_{\gamma}u\cdot\underbrace{\left( (\gamma^{ij}-\delta_{ij})u_k- \gamma^{ij} D\p(u)\left( g_k  \right)\right)}_{=:w_k^{ij}} \xi^k_{ij}\f \mu_{\gamma}=2\int_{B^+}\Delta_{\gamma}u\cdot w_k^{ij} \xi^k_{ij}\f \mu_{\gamma}\\
&=\dfrac{2}{\tau^2}\int_{B^+}\dfrac{\Delta_{\gamma}u\cdot w^{ij}_kx^ix^jx^k}{\vert x\vert^2} \si''\f \mu_{\gamma}\\
&\quad+\dfrac{2}{\tau}\int_{B^+}\Delta_{\gamma}u\cdot \left(w^{ik}_kx^i+w^{kj}_kx^j+w^{jj}_kx^k\right) \si'\cdot \dfrac{1}{\vert x\vert}\f \mu_{\gamma}\\
&\quad-\dfrac{2}{\tau}\int_{B^+}\dfrac{\Delta_{\gamma}u\cdot w_k^{ij}x^ix^jx^k}{\vert x\vert^2}\si'\cdot \dfrac{1}{\vert x\vert}\f \mu_{\gamma}
\numberthis\label{133}.
\end{align*}
We multiply \eqref{133} with $\e^{\chi\tau}\tau^{3-m}$ and integrate over $[\rho,r]$, i.e.,
\begin{align*}
\int_{\rho}^r\e^{\chi\tau}\tau^{3-m}VIII\f\tau&=2\int_{\rho}^r\e^{\chi\tau}\tau^{1-m}\int_{B^+}\dfrac{\Delta_{\gamma}u\cdot w^{ij}_kx^ix^jx^k}{\vert x\vert^2} \si''\f \mu_{\gamma}\f\tau\\
&\quad+2\int_{\rho}^r\e^{\chi\tau}\tau^{2-m}\int_{B^+}\Delta_{\gamma}u\cdot \left(w^{ik}_kx^i+w^{kj}_kx^j+w^{jj}_kx^k\right) \si'\cdot \dfrac{1}{\vert x\vert}\f \mu_{\gamma}\f\tau\\
&\quad-2\int_{\rho}^r\e^{\chi\tau}\tau^{2-m}\int_{B^+}\dfrac{\Delta_{\gamma}u\cdot w_k^{ij}x^ix^jx^k}{\vert x\vert^2}\si'\cdot \dfrac{1}{\vert x\vert}\f \mu_{\gamma}\f\tau\numberthis\label{134}.
\end{align*}
We reform the first integral in \eqref{134} with Fubini and integration by parts as follows,
\begin{align*}
&2\int_{\rho}^r\e^{\chi\tau}\tau^{1-m}\int_{B^+}\dfrac{\Delta_{\gamma}u\cdot w^{ij}_kx^ix^jx^k}{\vert x\vert^2} \si''\f \mu_{\gamma}\f\tau\\
&=-2\int_{B^+}\dfrac{\Delta_{\gamma}u\cdot w^{ij}_kx^ix^jx^k}{\vert x\vert^2} \left(\e^{\chi r}r^{3-m}\psi'_{\nu,r}-\e^{\chi \rho}\rho^{3-m}\psi'_{\nu,\rho}\right)\dfrac{1}{\vert x\vert}\f \mu_{\gamma}\\
&\quad+2(3-m)\int_{\rho}^r\e^{\chi\tau}\tau^{2-m}\int_{B^+}\dfrac{\Delta_{\gamma}u\cdot w^{ij}_kx^ix^jx^k}{\vert x\vert^2} \si'\cdot\dfrac{1}{\vert x\vert}\f \mu_{\gamma}\f\tau\\
&\quad+2\chi \int_{\rho}^r\e^{\chi\tau}\tau^{3-m}\int_{B^+}\dfrac{\Delta_{\gamma}u\cdot w^{ij}_kx^ix^jx^k}{\vert x\vert^2} \si'\cdot\dfrac{1}{\vert x\vert}\f \mu_{\gamma}\f\tau\numberthis\label{135}.
\end{align*}
Putting \eqref{135} into \eqref{134} yields
\begin{align*}
&\int_{\rho}^r\e^{\chi\tau}\tau^{3-m}VIII\f\tau\\
&=-2\int_{B^+}\dfrac{\Delta_{\gamma}u\cdot w^{ij}_kx^ix^jx^k}{\vert x\vert^2} \left(\e^{\chi r}r^{3-m}\psi'_{\nu,r}-\e^{\chi \rho}\rho^{3-m}\psi'_{\nu,\rho}\right)\cdot\dfrac{1}{\vert x\vert}\f \mu_{\gamma}\\
&\quad+2(2-m)\int_{\rho}^r\e^{\chi\tau}\tau^{2-m}\int_{B^+}\dfrac{\Delta_{\gamma}u\cdot w^{ij}_kx^ix^jx^k}{\vert x\vert^2} \si'\cdot\dfrac{1}{\vert x\vert}\f \mu_{\gamma}\f\tau\\
&\quad+2\chi \int_{\rho}^r\e^{\chi\tau}\tau^{3-m}\int_{B^+}\dfrac{\Delta_{\gamma}u\cdot w^{ij}_kx^ix^jx^k}{\vert x\vert^2} \si'\cdot\dfrac{1}{\vert x\vert}\f \mu_{\gamma}\f\tau\\
&\quad+2\int_{\rho}^r\e^{\chi\tau}\tau^{2-m}\int_{B^+}\Delta_{\gamma}u\cdot \left(w^{ik}_kx^i+w^{kj}_kx^j+w^{jj}_kx^k\right) \si'\cdot \dfrac{1}{\vert x\vert}\f \mu_{\gamma}\f\tau\numberthis\label{136}.
\end{align*}
We obtain by Lebesgue's differentiation theorem and Lemma 2 in the appendix of \cite{1} as $\nu\rightarrow\infty$
\begin{align*}
&\lim_{\nu\rightarrow\infty}\int_{\rho}^r\e^{\chi\tau}\tau^{3-m}VIII\f\tau\\
&=-2\int_{S^+_r\setminus S^+_{\rho}}\e^{\chi\vert x\vert}\dfrac{\Delta_{\gamma}u\cdot w^{ij}_kx^ix^jx^k}{\vert x\vert^{m-1}} \sqrt{\gamma}\f \h^{m-1}+2(2-m)\int_{ B^+_r\setminus  B^+_{\rho}}\e^{\chi\vert x\vert}\dfrac{\Delta_{\gamma}u\cdot w^{ij}_kx^ix^jx^k}{\vert x\vert^{m}} \f \mu_{\gamma}\\
&\quad+2\chi \int_{ B^+_r\setminus  B^+_{\rho}}\e^{\chi\vert x\vert}\dfrac{\Delta_{\gamma}u\cdot w^{ij}_kx^ix^jx^k}{\vert x\vert^{m-1}} \f \mu_{\gamma}+2\int_{ B^+_r\setminus  B^+_{\rho}}\e^{\chi\vert x\vert}\dfrac{\Delta_{\gamma}u\cdot \left(w^{ik}_kx^i+w^{kj}_kx^j+w^{jj}_kx^k\right)}{\vert x\vert^{m-2}}\f \mu_{\gamma}\numberthis\label{137}
\end{align*}
for a.e. $0<\rho<r<R$. Thus, the following estimate holds
\begin{align*}
RHS\eqref{137}&\leq 2\int_{S^+_r\cup S_{\rho}^+}\e^{\chi\vert x\vert}\dfrac{\vert \Delta_{\gamma}u\vert \vert w^{ij}_k\vert}{\vert x\vert^{m-4}} \sqrt{\gamma}\f \h^{m-1}+2(m+1)\int_{ B^+_r\setminus B_{\rho}^+}\e^{\chi\vert x\vert}\dfrac{\vert \Delta_{\gamma}u\vert \vert  w^{ij}_k\vert }{\vert x\vert^{m-3}}\f \mu_{\gamma}\\
&\quad+2\chi\int_{ B^+_r\setminus B_{\rho}^+}\e^{\chi\vert x\vert}\dfrac{\vert \Delta_{\gamma}u\vert \vert  w^{ij}_k\vert }{\vert x\vert^{m-4}}\f \mu_{\gamma} \numberthis\label{140}.
\end{align*}
Since $\vert w^{ij}_k\vert\leq H\vert Du\vert \vert x\vert +C_{46}$ where $C_{46}=C_{46}(\Vert Dg\Vert_{\infty}):=G\cdot \sup_{B^+}\vert D\p\circ u\vert \Vert Dg\Vert_{\infty}$ and with $\vert \Delta_{\gamma}u\vert\leq G\vert D^2u\vert +C_7\vert Du\vert$ we obtain from \eqref{140}
\begin{align*}
&\lim_{\nu\rightarrow\infty}\int_{\rho}^r\e^{\chi\tau}\tau^{3-m}VII\f\tau\leq RHS\eqref{140}\\
&\leq C_{47}\int_{S^+_r\cup S_{\rho}^+}\e^{\chi\vert x\vert}\dfrac{\vert D^2u\vert^2 }{\vert x\vert^{m-6}} \sqrt{\gamma}\f \h^{m-1}+C_{48}\int_{S^+_r\cup S_{\rho}^+}\e^{\chi\vert x\vert}\dfrac{\vert Du\vert^2 }{\vert x\vert^{m-4}} \sqrt{\gamma}\f \h^{m-1}\\
&\quad+C_{49}\int_{ B^+_r\setminus B_{\rho}^+}\e^{\chi\vert x\vert}\dfrac{\vert D^2u\vert^2  }{\vert x\vert^{m-5}}\f \mu_{\gamma}+C_{50}\int_{ B^+_r\setminus B_{\rho}^+}\e^{\chi\vert x\vert}\dfrac{\vert Du\vert^2 }{\vert x\vert^{m-3}}\f \mu_{\gamma}\\
&\quad+C_{45}\cdot (r+\rho)+C_{51}\cdot (r-\rho)\numberthis\label{143}
\end{align*}
where $C_{45}=C_{45}(\Vert Dg\Vert_{C^2},\Vert u\Vert_{L^{\infty}(B_1^+)})$,  $C_{47}=C_{47}(\Vert Dg\Vert_{\infty})$, $C_{48}=C_{48}(\Vert Dg\Vert_{\infty})$, $C_{49}=C_{49}( \Vert Dg\Vert_{C^2},\Vert u\Vert_{L^{\infty}(B_1^+)})$, $C_{50}=C_{50}(\Vert Dg\Vert_{C^2},\Vert u\Vert_{L^{\infty}(B_1^+)})$ and
$C_{51}=C_{51}(\Vert Dg\Vert_{C^2},\Vert u\Vert_{L^{\infty}(B_1^+)})$. Altogether, it holds because of \eqref{131} and \eqref{143} (cf. \eqref{83})
\begin{align*}
&\lim_{\nu\rightarrow\infty}\int_{r}^{\rho}\e^{\chi \tau}\tau^{3-m}LHS\eqref{58}\f\tau\\
&\leq \int_{\rho}^r\e^{\chi\tau}I'(\tau)\f\tau+2C_7\int_{\rho}^r\e^{\chi\tau}I(\tau)\f\tau+C_{53}\cdot (r-\rho)+C_{45}\cdot (r+\rho)\\
&\quad+C_{42}\int_{\rho}^r\e^{\chi\tau}\tau^{4-m}\int_{B^+_{\tau}} \vert D^2u\vert^2 \f\mu_{\gamma}\f\tau\\
&\quad-2\int_{S^+_r\setminus S^+_{\rho}}\e^{\chi\vert x\vert}\left(-\dfrac{ u_{i} u_{ik}x^k}{\vert x\vert^{m-3}}+ 2\dfrac{ (u_{i}x^i)^2}{\vert x\vert^{m-1}}-2\dfrac{ (u_{i})^2}{\vert x\vert^{m-3}}\right)\sqrt{\gamma}\f\h^{m-1}\\
&\quad-4\int_{B^+_r\setminus B^+_{\rho}}\e^{\chi\vert x\vert}\left( \dfrac{ (u_i+u_{ij}x^j)^2}{\vert x\vert^{m-2}}+\dfrac{(m-2)(u_{i}x^i)^2}{\vert x\vert^{m}}\right)\f\mu_{\gamma}\\
&\quad+2\chi\int_{ B^+_{r}\setminus  B^+_{\rho}}\e^{\chi \vert x\vert}\left(-\dfrac{ u_{i} u_{ik}x^k}{\vert x\vert^{m-3}}+2\dfrac{ (u_{i} x^i)^2}{\vert x\vert^{m-1}}-2\dfrac{ (u_{i})^2}{\vert x\vert^{m-3}}\right)\f\mu_{\gamma}\\
&\quad +C_{54}\int_{ B^+_{r}\setminus B_{\rho}^+}\e^{\chi\vert x\vert}\dfrac{\vert D^2 u\vert^2  }{\vert x\vert^{m-5}}\f\mu_{\gamma}+C_{55}\int_{ B^+_{r}\setminus B_{\rho}^+}\e^{\chi\vert x\vert}\dfrac{\vert Du\vert^2  }{\vert x\vert^{m-3}}\f\mu_{\gamma}\\
&\quad+C_{56}\cdot \int_{S_{r}^+\cup S_{\rho}^+}\e^{\chi \vert x\vert}\dfrac{\vert D^2u\vert^2  }{\vert x\vert^{m-6}} \sqrt{\gamma}\f\h^{m-1}+C_{57}\cdot \int_{S_{r}^+\cup S_{\rho}^+}\e^{\chi \vert x\vert}\dfrac{\vert Du\vert^2}{\vert x\vert^{m-4}}\sqrt{\gamma} \f\h^{m-1}\numberthis\label{144}.
\end{align*}
where $C_{53}=C_{53}(\Vert Dg\Vert_{C^2},\Vert u\Vert_{L^{\infty}(B_1^+)}):=C_{39}+C_{51}$, $C_{54}=C_{54}(\Vert Dg\Vert_{C^2},\Vert u\Vert_{L^{\infty}(B_1^+)}):=C_{40}+C_{49}$, $C_{55}=C_{55}(\Vert Dg\Vert_{C^2},\Vert u\Vert_{L^{\infty}(B_1^+)}):=C_{41}+C_{50}$, $C_{56}=C_{56}(\Vert Dg\Vert_{\infty}):=H+C_{47}$ and $C_{57}=C_{57}(\Vert Dg\Vert_{\infty}):=C_{48}+2C_7+H$. Now, from \eqref{82} we infer the chain of inequalities
\begin{align*}
RHS\eqref{144}\geq \lim_{\nu\rightarrow\infty}\int_{r}^{\rho}\e^{\chi \tau}\tau^{3-m}LHS\eqref{58}\f\tau\geq -RHS\eqref{81}\numberthis\label{145},
\end{align*}
i.e.
\begin{align*}
&\int_{\rho}^r\e^{\chi\tau}I'(\tau)\f\tau+2C_7\int_{\rho}^r\e^{\chi\tau}I(\tau)\f\tau+C_{42}\int_{\rho}^r\e^{\chi\tau}J(\tau)\f\tau+2C_{45}\cdot r\\
&\quad-2\int_{S^+_r\setminus S^+_{\rho}}\e^{\chi\vert x\vert}\left(-\dfrac{ u_{i} u_{ik}x^k}{\vert x\vert^{m-3}}+ 2\dfrac{ (u_{i}x^i)^2}{\vert x\vert^{m-1}}-2\dfrac{ (u_{i})^2}{\vert x\vert^{m-3}}\right)\sqrt{\gamma}\f\h^{m-1}\\
&\quad-4\int_{B^+_r\setminus B^+_{\rho}}\e^{\chi\vert x\vert}\left( \dfrac{ (u_i+u_{ij}x^j)^2}{\vert x\vert^{m-2}}+\dfrac{(m-2)(u_{i}x^i)^2}{\vert x\vert^{m}}\right)\f\mu_{\gamma}\\
&\quad+2\chi\int_{ B^+_{r}\setminus  B^+_{\rho}}\e^{\chi \vert x\vert}\left(-\dfrac{ u_{i} u_{ik}x^k}{\vert x\vert^{m-3}}+2\dfrac{ (u_{i} x^i)^2}{\vert x\vert^{m-1}}-2\dfrac{ (u_{i})^2}{\vert x\vert^{m-3}}\right)\f\mu_{\gamma}\\
&\quad +C_{54}\int_{ B^+_{r}\setminus B_{\rho}^+}\e^{\chi \vert x\vert}\dfrac{\vert D^2 u\vert^2  }{\vert x\vert^{m-5}}\f\mu_{\gamma}+C_{55}\int_{ B^+_{r}\setminus B_{\rho}^+}\e^{\chi \vert x\vert}\dfrac{\vert Du\vert^2  }{\vert x\vert^{m-3}}\f\mu_{\gamma}\\
&\quad+C_{56}\cdot \int_{S_{r}^+\cup S_{\rho}^+}\e^{\chi \vert x\vert}\dfrac{\vert D^2u\vert^2  }{\vert x\vert^{m-6}} \sqrt{\gamma}\f\h^{m-1}+C_{57}\cdot \int_{S_{r}^+\cup S_{\rho}^+}\e^{\chi \vert x\vert}\dfrac{\vert Du\vert^2}{\vert x\vert^{m-4}}\sqrt{\gamma} \f\h^{m-1}\\
&\geq -C_{22} \int_{\rho}^r\e^{\chi\tau}I(\tau)\f\tau-C_{27}\int_{\rho}^r\e^{\chi\tau}J(\tau)\f\tau-C_{30}\cdot (r-\rho)\\
&\quad-C_{20}\int_{B^+_r\setminus B_{\rho}^+}\e^{\chi \vert x\vert}\dfrac{\vert D^2u\vert^2}{\vert x\vert^{m-5}}  \f \mu_{\gamma}-C_{31}\int_{B^+_{r}\setminus B_{\rho}^+}\e^{\chi \vert x\vert}\dfrac{\vert Du\vert^2}{\vert x\vert^{m-3}} \f \mu_{\gamma}\numberthis\label{146}
\end{align*}
where we set $J(\tau):=\displaystyle\tau^{4-m}\int_{B^+_{\tau}} \vert D^2u\vert^2 \f\mu_{\gamma}$ and used $C_{53}\cdot (r-\rho)+C_{45}\cdot (r+\rho)\leq 2C_{45}\cdot r$, since $\rho < r$. This inequality can be rewritten to
\begin{align*}
&\int_{\rho}^r\e^{\chi\tau}I'(\tau)\f\tau+(2C_7+C_{22})\int_{\rho}^r\e^{\chi\tau}I(\tau)\f\tau+(C_{30}+C_{32}+2C_{45})\cdot r-(C_{30}+C_{32})\cdot \rho\\
&\quad+(C_{42}+C_{27})\int_{\rho}^r\e^{\chi\tau}J(\tau)\f\tau\\
&\quad +(C_{20}+C_{54})\cdot \int_{ B^+_{r}\setminus B_{\rho}^+}\e^{\chi \vert x\vert}\dfrac{\vert D^2 u\vert^2  }{\vert x\vert^{m-5}}\f\mu_{\gamma}+(C_{31}+C_{55})\cdot \int_{ B^+_{r}\setminus B_{\rho}^+}\e^{\chi \vert x\vert}\dfrac{\vert Du\vert^2  }{\vert x\vert^{m-3}}\f\mu_{\gamma}\\
&\quad+C_{56}\cdot \int_{S_{r}^+\cup S_{\rho}^+}\e^{\chi \vert x\vert}\dfrac{\vert D^2u\vert^2  }{\vert x\vert^{m-6}} \sqrt{\gamma}\f\h^{m-1}+C_{57}\cdot \int_{S_{r}^+\cup S_{\rho}^+}\e^{\chi \vert x\vert}\dfrac{\vert Du\vert^2}{\vert x\vert^{m-4}}\sqrt{\gamma} \f\h^{m-1}\\
&\geq 4\int_{B^+_r\setminus B^+_{\rho}}\e^{\chi\vert x\vert}\left( \dfrac{ (u_i+u_{ij}x^j)^2}{\vert x\vert^{m-2}}+\dfrac{(m-2)(u_{i}x^i)^2}{\vert x\vert^{m}}\right)\f\mu_{\gamma}\\
&\quad+2\int_{S^+_r\setminus S^+_{\rho}}\e^{\chi\vert x\vert}\left(-\dfrac{ u_{i} u_{ik}x^k}{\vert x\vert^{m-3}}+ 2\dfrac{ (u_{i}x^i)^2}{\vert x\vert^{m-1}}-2\dfrac{ (u_{i})^2}{\vert x\vert^{m-3}}\right)\sqrt{\gamma}\f\h^{m-1}\\
&\quad-2\chi\int_{ B^+_{r}\setminus  B^+_{\rho}}\e^{\chi \vert x\vert}\left(-\dfrac{ u_{i} u_{ik}x^k}{\vert x\vert^{m-3}}+2\dfrac{ (u_{i} x^i)^2}{\vert x\vert^{m-1}}-2\dfrac{ (u_{i})^2}{\vert x\vert^{m-3}}\right)\f\mu_{\gamma}\numberthis\label{147}.
\end{align*}

\textit{Step 7}. Next, we recast $J(\tau)=\displaystyle\tau^{4-m}\int_{B^+_{\tau}} \vert D^2u\vert^2 \f\mu_{\gamma}$ with the aid of integration by parts. For that purpose we approximate $u$ by $u^{\varepsilon}$ as in \Cref{lem2.6}. It holds
\begin{align*}
&\tau^{4-m}\int_{B^+_{\tau}}\vert D^2u^{\varepsilon}\vert^2 \f\mu_{\gamma}\\
&=\tau^{4-m}\int_{\partial B^+_{\tau}} u^{\varepsilon}_i u^{\varepsilon}_{ij}\nu^j  \sqrt{\gamma}\f\h^{m-1}-\tau^{4-m}\int_{B^+_{\tau}}u^{\varepsilon}_iu^{\varepsilon}_{ijj} \f\mu_{\gamma}-\tau^{4-m}\int_{B^+_{\tau}} u^{\varepsilon}_i u^{\varepsilon}_{ij}\left(\sqrt{\gamma}\right)_j \f\mathcal{L}^m\\
&=\tau^{4-m}\int_{\partial B^+_{\tau}} \left(u^{\varepsilon}_i u^{\varepsilon}_{ij}\nu^j-u^{\varepsilon}_i\nu^iu^{\varepsilon}_{jj}\right) \sqrt{\gamma}\f\h^{m-1}\\
&\quad+\tau^{4-m}\int_{B^+_{\tau}}\left(u^{\varepsilon}_i\dfrac{\left(\sqrt{\gamma}\right)_i}{\sqrt{\gamma}} u^{\varepsilon}_{jj}- u^{\varepsilon}_i u^{\varepsilon}_{ij}\dfrac{\left(\sqrt{\gamma}\right)_j}{\sqrt{\gamma}}\right) \f\mu_{\gamma}+\tau^{4-m}\int_{B^+_{\tau}}\vert \Delta u^{\varepsilon}\vert^2 \f\mu_{\gamma}\\
&=\int_{S^+_{\tau}} \left(\dfrac{u^{\varepsilon}_i u^{\varepsilon}_{ij}x^j}{\vert x\vert^{m-3}}-\dfrac{u^{\varepsilon}_ix^iu^{\varepsilon}_{jj}}{\vert x\vert^{m-3}}\right) \sqrt{\gamma}\f\h^{m-1}-\tau^{4-m}\int_{T_{\tau}} \left(g_i g_{im}-g_mg_{jj}\right) \sqrt{\gamma}\f\h^{m-1}\\
&\quad+\tau^{4-m}\int_{B^+_{\tau}}\left(u^{\varepsilon}_i\dfrac{\left(\sqrt{\gamma}\right)_i}{\sqrt{\gamma}} u^{\varepsilon}_{jj}- u^{\varepsilon}_i u^{\varepsilon}_{ij}\dfrac{\left(\sqrt{\gamma}\right)_j}{\sqrt{\gamma}}\right) \f\mu_{\gamma}+\tau^{4-m}\int_{B^+_{\tau}}\vert \Delta u^{\varepsilon}\vert^2 \f\mu_{\gamma}\\
&\leq \int_{S^+_{\tau}} \left(\dfrac{u^{\varepsilon}_i u^{\varepsilon}_{ij}x^j}{\vert x\vert^{m-3}}-\dfrac{u^{\varepsilon}_ix^iu^{\varepsilon}_{jj}}{\vert x\vert^{m-3}}\right) \sqrt{\gamma}\f\h^{m-1} +C_{58}R^3\\
&\quad+\tau^{4-m}\int_{B^+_{\tau}}\left(u^{\varepsilon}_{jj}u^{\varepsilon}_i\Gamma^l_{il} - u^{\varepsilon}_i u^{\varepsilon}_{ij}\Gamma^l_{jl}\right) \f\mu_{\gamma}+\tau^{4-m}\int_{B^+_{\tau}}\vert \Delta u^{\varepsilon}\vert^2 \f\mu_{\gamma} \numberthis\label{148}
\end{align*}
where $C_{58}=C_{58}(\Vert Dg\Vert_{C^1}):=2\Vert Dg\Vert_{\infty}\Vert D^2g\Vert_{\infty}G^{m/2}\h^{m-1}(T_1)$.  As $\varepsilon\searrow0$, we obtain 
\begin{align*}
\tau^{4-m}\int_{B^+_{\tau}}\vert D^2u\vert^2 \f\mu_{\gamma}&\leq \int_{S^+_{\tau}} \left(\dfrac{u_i u_{ij}x^j}{\vert x\vert^{m-3}}-\dfrac{u_ix^iu_{jj}}{\vert x\vert^{m-3}}\right) \sqrt{\gamma}\f\h^{m-1} +C_{59}\\
&\quad+\tau^{4-m}\int_{B^+_{\tau}}\left(u_{jj}u_i\Gamma^l_{il} - u_i u_{ij}\Gamma^l_{jl}\right) \f\mu_{\gamma}+\tau^{4-m}\int_{B^+_{\tau}}\vert \Delta u\vert^2 \f\mu_{\gamma} \numberthis\label{149}
\end{align*}
for a.e. $\tau\in (0,R)$ where $C_{59}(\Vert Dg\Vert_{C^1},R):=C_{58}R^3$. Notice that $\Delta u=\Delta_{\gamma}u-(\gamma^{ij}-\delta^{ij})u_{ij}+\gamma^{ij}\Gamma^k_{ij}u_k$ holds. Thus, $\vert \Delta u\vert^2\leq 3\vert\Delta_{\gamma}u\vert^2 + 3H^2\vert x\vert^2 \vert D^2u\vert^2 +3C_7^2\vert Du\vert^2$. Moreover, we estimate for the second to last integral in \eqref{149} with $2\vert D^2u\vert \vert Du\vert\leq \tau \vert D^2u\vert^2 +\dfrac{1}{\tau}\vert Du\vert^2$ as follows,
\begin{align*}
&\tau^{4-m}\int_{B^+_{\tau}}\left(u_{jj}u_i\Gamma^l_{il} - u_i u_{ij}\Gamma^l_{jl}\right) \f\mu_{\gamma}\\
&\leq C_5\tau^{5-m}\int_{B^+_{\tau}}\vert D^2u\vert^2 \f\mu_{\gamma}+C_5\tau^{3-m}\int_{B^+_{\tau}} \vert Du\vert^2 \f\mu_{\gamma}\numberthis\label{152}.
\end{align*}
Thus, it holds
\begin{align*}
&\int_{\rho}^r\e^{\chi\tau}\tau^{4-m}\int_{B^+_{\tau}} \vert D^2u\vert^2 \f\mu_{\gamma}\f\tau\\
&\leq \int_{B^+_{r}\setminus B_{\rho}^+}\e^{\chi\vert x\vert} \left(\dfrac{u_i u_{ij}x^j}{\vert x\vert^{m-3}}-\dfrac{u_ix^iu_{jj}}{\vert x\vert^{m-3}}\right) \f \mu_{\gamma}+ 3\int_{\rho}^r\e^{\chi\tau}I(\tau)\f\tau +C_{59}\e^{\chi R}\cdot (r-\rho)\\
&\quad+C_{60}\int_{\rho}^r\e^{\chi\tau}\tau^{5-m}\int_{B^+_{\tau}}\vert D^2u\vert^2 \f\mu_{\gamma}\f\tau+ C_{61}\int_{\rho}^r\e^{\chi\tau}\tau^{3-m}\int_{B^+_{\tau}} \vert Du\vert^2 \f\mu_{\gamma}\f\tau\numberthis\label{153}
\end{align*}
where $C_{60}:=C_5+3H^2$ and $C_{61}:=C_5+3C_7^2$.\\

\textit{Step 8}. We put $\chi=\chi(\Vert Dg\Vert_{C^2},\Vert u\Vert_{L^{\infty}}):=2C_7+C_{22}+3C_{27}+3C_{42}$. Then, we obtain from \eqref{147} with the aid of \eqref{153} the inequality
\begin{align*}
&\int_{\rho}^r\e^{\chi\tau}I'(\tau)\f\tau+\chi\int_{\rho}^r\e^{\chi\tau}I(\tau)\f\tau+\tilde{C}_{62}\cdot r-\tilde{C}_{63}\cdot \rho\\
&\quad+(C_{27}+C_{42})\left(C_{60}\int_{\rho}^r\e^{\chi\tau}\tau^{5-m}\int_{B^+_{\tau}}\vert D^2u\vert^2 \f\mu_{\gamma}\f\tau+ C_{61}\int_{\rho}^r\e^{\chi\tau}\tau^{3-m}\int_{B^+_{\tau}} \vert Du\vert^2 \f\mu_{\gamma}\f\tau\right)\\
&\quad +\tilde{C}_{60}\cdot\int_{ B^+_{r}\setminus B_{\rho}^+}\e^{\chi \vert x\vert}\dfrac{\vert D^2 u\vert^2  }{\vert x\vert^{m-5}}\f\mu_{\gamma}+\tilde{C}_{61}\cdot \int_{ B^+_{r}\setminus B_{\rho}^+}\e^{\chi \vert x\vert}\dfrac{\vert Du\vert^2  }{\vert x\vert^{m-3}}\f\mu_{\gamma}\\
&\quad+C_{56}\cdot\int_{S_{r}^+\cup S_{\rho}^+}\e^{\chi \vert x\vert}\dfrac{\vert D^2u\vert^2  }{\vert x\vert^{m-6}} \sqrt{\gamma}\f\h^{m-1}+C_{57}\cdot \int_{S_{r}^+\cup S_{\rho}^+}\e^{\chi \vert x\vert}\dfrac{\vert Du\vert^2}{\vert x\vert^{m-4}}\sqrt{\gamma} \f\h^{m-1}\\
&\geq 4\int_{B^+_r\setminus B^+_{\rho}}\e^{\chi\vert x\vert}\left( \dfrac{ (u_i+u_{ij}x^j)^2}{\vert x\vert^{m-2}}+\dfrac{(m-2)(u_{i}x^i)^2}{\vert x\vert^{m}}\right)\f\mu_{\gamma}\\
&\quad+2\int_{S^+_r\setminus S^+_{\rho}}\e^{\chi\vert x\vert}\left(-\dfrac{ u_{i} u_{ik}x^k}{\vert x\vert^{m-3}}+ 2\dfrac{ (u_{i}x^i)^2}{\vert x\vert^{m-1}}-2\dfrac{ (u_{i})^2}{\vert x\vert^{m-3}}\right)\sqrt{\gamma}\f\h^{m-1}\\
&\quad+2\chi\int_{ B^+_{r}\setminus  B^+_{\rho}}\e^{\chi \vert x\vert}\left(\dfrac{ u_{i} u_{ik}x^k}{\vert x\vert^{m-3}}-2\dfrac{ (u_{i} x^i)^2}{\vert x\vert^{m-1}}+2\dfrac{ (u_{i})^2}{\vert x\vert^{m-3}}\right)\f\mu_{\gamma}\\
&\quad-(C_{27}+C_{42})\int_{B^+_{r}\setminus B_{\rho}^+}\e^{\chi\vert x\vert} \left(\dfrac{u_i u_{ij}x^j}{\vert x\vert^{m-3}}-\dfrac{u_ix^iu_{jj}}{\vert x\vert^{m-3}}\right) \f \mu_{\gamma}\numberthis\label{154}
\end{align*}
where $\tilde{C}_{60}=\tilde{C}_{60}(\Vert Dg\Vert_{C^2},\Vert u\Vert_{L^{\infty}(B_1^+)}):=C_{20}+C_{54}$, $\tilde{C}_{61}=\tilde{C}_{61}(\Vert Dg\Vert_{C^2},\Vert u\Vert_{L^{\infty}(B_1^+)}):=C_{31}+C_{55}$, $\tilde{C}_{62}=\tilde{C}_{62}(\Vert Dg\Vert_{C^2},\Vert u\Vert_{L^{\infty}(B_1^+)}):=C_{30}+C_{32}+2C_{45}+(C_{27}+C_{42})C_{59}\e^{\chi R}$, $\tilde{C}_{63}=\tilde{C}_{63}(\Vert Dg\Vert_{C^2},\Vert u\Vert_{L^{\infty}(B_1^+)}):=C_{30}+C_{32}+(C_{27}+C_{42})C_{59}\e^{\chi R}$ with $\tilde{C}_{62}\geq \tilde{C}_{63}$. Beyond that, it holds 
\begin{align*}
-2\dfrac{ (u_{i} x^i)^2}{\vert x\vert^{m-1}}+2\dfrac{ (u_{i})^2}{\vert x\vert^{m-3}}\geq -2\dfrac{ \vert Du\vert^2 \vert x \vert^2}{\vert x\vert^{m-1}}+2\dfrac{ \vert Du\vert^2}{\vert x\vert^{m-3}}=0\numberthis\label{155}.
\end{align*}
Hence, applying Young's inequality we obtain
\begin{align*}
&2\chi\int_{ B^+_{r}\setminus  B^+_{\rho}}\e^{\chi \vert x\vert}\left(\dfrac{ u_{i} u_{ik}x^k}{\vert x\vert^{m-3}}-2\dfrac{ (u_{i} x^i)^2}{\vert x\vert^{m-1}}+2\dfrac{ (u_{i})^2}{\vert x\vert^{m-3}}\right)\f\mu_{\gamma}\\
&\geq -\chi \int_{ B^+_{r}\setminus  B^+_{\rho}}\e^{\chi \vert x\vert}\dfrac{ \vert D^2 u\vert^2   }{\vert x\vert^{m-5}}\f\mu_{\gamma}-\chi \int_{ B^+_{r}\setminus  B^+_{\rho}}\e^{\chi \vert x\vert}\dfrac{ \vert D u\vert^2   }{\vert x\vert^{m-3}}\f\mu_{\gamma}\numberthis\label{156}.
\end{align*}
Analogously, we obtain for the last integral on the left-hand side of \eqref{154} the estimate
\begin{align*}
&-(C_{27}+C_{42})\int_{B^+_{r}\setminus B_{\rho}^+}\e^{\chi\vert x\vert} \left(\dfrac{u_i u_{ij}x^j}{\vert x\vert^{m-3}}-\dfrac{u_ix^iu_{jj}}{\vert x\vert^{m-3}}\right) \f \mu_{\gamma}\\
&\geq -(C_{27}+C_{42}) \int_{ B^+_{r}\setminus B_{\rho}^+}\e^{\chi \vert x\vert}\dfrac{ \vert D^2 u\vert^2   }{\vert x\vert^{m-5}}\f\mu_{\gamma}-(C_{27}+C_{42}) \int_{ B^+_{r}\setminus B_{\rho}^+}\e^{\chi \vert x\vert}\dfrac{ \vert D u\vert^2   }{\vert x\vert^{m-3}}\f\mu_{\gamma}\numberthis\label{157}.
\end{align*}
Moreover, observe that it holds
\begin{align*}
\int_{\rho}^r\e^{\chi\tau}I'(\tau)\f\tau+\chi\int_{\rho}^r\e^{\chi\tau}I(\tau)\f\tau&=\int_{\rho}^r\dfrac{\f}{\f \tau}\left(\e^{\chi\tau}I(\tau)\right)\f\tau=\e^{\chi r}I(r)-\e^{\chi\rho}I(\rho)\numberthis\label{159}.
\end{align*}
Because of \eqref{156}, \eqref{157}, \eqref{159} and $\rho < r$ we obtain from \eqref{154} the inequality
\begin{align*}
&\e^{\chi r}r^{4-m}\int_{B_r^+}\vert \Delta_{\gamma}u\vert^2 \f\mu_{\gamma}-\e^{\chi\rho}\rho^{4-m}\int_{B_{\rho}^+}\vert \Delta_{\gamma}u\vert^2 \f\mu_{\gamma}+\mathsf{C}_1\cdot r\\
&\quad+\mathsf{C}_2\int_{\rho}^r\e^{\chi\tau}\tau^{5-m}\int_{B^+_{\tau}}\vert D^2u\vert^2 \f\mu_{\gamma}\f\tau+ \mathsf{C}_2\int_{\rho}^r\e^{\chi\tau}\tau^{3-m}\int_{B^+_{\tau}} \vert Du\vert^2 \f\mu_{\gamma}\f\tau\\
&\quad +\mathsf{C}_3\int_{ B^+_{r}\setminus B_{\rho}^+}\e^{\chi \vert x\vert}\dfrac{\vert D^2 u\vert^2  }{\vert x\vert^{m-5}}\f\mu_{\gamma}+\mathsf{C}_4\cdot \int_{ B^+_{r}\setminus B_{\rho}^+}\e^{\chi \vert x\vert}\dfrac{\vert Du\vert^2  }{\vert x\vert^{m-3}}\f\mu_{\gamma}\\
&\quad+\mathsf{C}_6\int_{S_{r}^+\cup S_{\rho}^+}\e^{\chi \vert x\vert}\dfrac{\vert D^2u\vert^2  }{\vert x\vert^{m-6}} \sqrt{\gamma}\f\h^{m-1}+\mathsf{C}_7 \int_{S_{r}^+\cup S_{\rho}^+}\e^{\chi \vert x\vert}\dfrac{\vert Du\vert^2}{\vert x\vert^{m-4}}\sqrt{\gamma} \f\h^{m-1}\\
&\geq 4\int_{B^+_r\setminus B^+_{\rho}}\e^{\chi\vert x\vert}\left( \dfrac{ (u_i+u_{ij}x^j)^2}{\vert x\vert^{m-2}}+\dfrac{(m-2)(u_{i}x^i)^2}{\vert x\vert^{m}}\right)\f\mu_{\gamma}\\
&\quad+2\int_{S^+_r\setminus S^+_{\rho}}\e^{\chi\vert x\vert}\left(-\dfrac{ u_{i} u_{ik}x^k}{\vert x\vert^{m-3}}+ 2\dfrac{ (u_{i}x^i)^2}{\vert x\vert^{m-1}}-2\dfrac{ (u_{i})^2}{\vert x\vert^{m-3}}\right)\sqrt{\gamma}\f\h^{m-1}\numberthis\label{161}
\end{align*}
where we have set $\mathsf{C}_1=\mathsf{C}_1(\Vert Dg\Vert_{C^2},\Vert u\Vert_{L^{\infty}}):=\tilde{C}_{62}$, $\mathsf{C}_2=\mathsf{C}_2(\Vert Dg\Vert_{C^1},\Vert u\Vert_{L^{\infty}}):=(C_{27}+C_{42})C_{60}$, $\mathsf{C}_3=\mathsf{C}_3(\Vert Dg\Vert_{C^1},\Vert u\Vert_{L^{\infty}}):=(C_{27}+C_{42})C_{61}$, $\mathsf{C}_4=\mathsf{C}_4(\Vert Dg\Vert_{C^2},\Vert u\Vert_{L^{\infty}}):=\tilde{C}_{60}+\chi+C_{27}+C_{42}$, $\mathsf{C}_5=\mathsf{C}_5(\Vert Dg\Vert_{C^2},\Vert u\Vert_{L^{\infty}}):=\tilde{C}_{61}+\chi+C_{27}+C_{42}$, $\mathsf{C}_6=\mathsf{C}_6(\Vert Dg\Vert_{\infty}):=C_{56}$ and  $\mathsf{C}_7=\mathsf{C}_7(\Vert Dg\Vert_{\infty}):=C_{57}$. So, we have the boundary monotonicity inequality \eqref{57} for $a=0$. This concludes the proof of \Cref{thm}.
\end{proof}

\textbf{Acknowledgements.} I would like to thank Prof. Dr. Christoph Scheven for his many helpful advices.

\end{document}